\documentclass[11pt]{article}
\usepackage{amsmath, upgreek}
\usepackage{amssymb}
\usepackage{color}
\usepackage{amscd}
\usepackage{xspace}
\usepackage{verbatim}
\usepackage{graphicx}
\newcommand{\mysection}[1]{
\section{#1}\setcounter{equation}{0}}
\title{\bf A priori estimates for elliptic equations
 with reaction terms involving  the function and its gradient}
\author{{\bf Marie-Fran\c{c}oise Bidaut-V\'eron\footnote{\noindent Laboratoire de Math\'{e}matiques et Physique Th\'{e}orique,
Universit\'e de Tours, 37200 Tours, France. E-mail: veronmf@univ-tours.fr},} \\{\bf Marta Garcia-Huidobro \footnote{\noindent
Departamento de Matematicas, Pontifica Universidad Catolica de Chile
Casilla 307, Correo 2, Santiago de Chile. E-mail: mgarcia@mat.puc.cl}}\\
 {\bf Laurent V\'eron \footnote{\noindent
Laboratoire de Math\'{e}matiques et Physique Th\'{e}orique, Universit\'e de Tours, 37200 Tours, France. E-mail: veronl@univ-tours.fr}}\\[2mm]
}

\date{}
\begin{document}
 \maketitle


\newcommand{\txt}[1]{\;\text{ #1 }\;}
\newcommand{\tbf}{\textbf}
\newcommand{\tit}{\textit}
\newcommand{\tsc}{\textsc}
\newcommand{\trm}{\textrm}
\newcommand{\mbf}{\mathbf}
\newcommand{\mrm}{\mathrm}
\newcommand{\bsym}{\boldsymbol}
\newcommand{\scs}{\scriptstyle}
\newcommand{\sss}{\scriptscriptstyle}
\newcommand{\txts}{\textstyle}
\newcommand{\dsps}{\displaystyle}
\newcommand{\fnz}{\footnotesize}
\newcommand{\scz}{\scriptsize}
\newcommand{\be}{\begin{equation}}
\newcommand{\bel}[1]{\begin{equation}\label{#1}}
\newcommand{\ee}{\end{equation}}
\newcommand{\eqnl}[2]{\begin{equation}\label{#1}{#2}\end{equation}}
\newcommand{\barr}{\begin{eqnarray}}
\newcommand{\earr}{\end{eqnarray}}
\newcommand{\bars}{\begin{eqnarray*}}
\newcommand{\ears}{\end{eqnarray*}}
\newcommand{\nnu}{\nonumber \\}
\newtheorem{subn}{\name}
\renewcommand{\thesubn}{}
\newcommand{\bsn}[1]{\def\name{#1}\begin{subn}}
\newcommand{\esn}{\end{subn}}
\newtheorem{sub}{\name}[section]
\newcommand{\dn}[1]{\def\name{#1}}   
\newcommand{\bs}{\begin{sub}}
\newcommand{\es}{\end{sub}}
\newcommand{\bsl}[1]{\begin{sub}\label{#1}}
\newcommand{\bth}[1]{\def\name{Theorem}
\begin{sub}\label{t:#1}}
\newcommand{\blemma}[1]{\def\name{Lemma}
\begin{sub}\label{l:#1}}
\newcommand{\bcor}[1]{\def\name{Corollary}
\begin{sub}\label{c:#1}}
\newcommand{\bdef}[1]{\def\name{Definition}
\begin{sub}\label{d:#1}}
\newcommand{\bprop}[1]{\def\name{Proposition}
\begin{sub}\label{p:#1}}

\newcommand{\R}{\eqref}
\newcommand{\rth}[1]{Theorem~\ref{t:#1}}
\newcommand{\rlemma}[1]{Lemma~\ref{l:#1}}
\newcommand{\rcor}[1]{Corollary~\ref{c:#1}}
\newcommand{\rdef}[1]{Definition~\ref{d:#1}}
\newcommand{\rprop}[1]{Proposition~\ref{p:#1}}
\newcommand{\BA}{\begin{array}}
\newcommand{\EA}{\end{array}}
\newcommand{\BAN}{\renewcommand{\arraystretch}{1.2}
\setlength{\arraycolsep}{2pt}\begin{array}}
\newcommand{\BAV}[2]{\renewcommand{\arraystretch}{#1}
\setlength{\arraycolsep}{#2}\begin{array}}
\newcommand{\BSA}{\begin{subarray}}
\newcommand{\ESA}{\end{subarray}}
\newcommand{\BAL}{\begin{aligned}}
\newcommand{\EAL}{\end{aligned}}
\newcommand{\BALG}{\begin{alignat}}
\newcommand{\EALG}{\end{alignat}}
\newcommand{\BALGN}{\begin{alignat*}}
\newcommand{\EALGN}{\end{alignat*}}
\newcommand{\note}[1]{\textit{#1.}\hspace{2mm}}
\newcommand{\Proof}{\note{Proof}}
\newcommand{\qeda}{\hspace{10mm}\hfill $\square$}
\newcommand{\qed}{\\
${}$ \hfill $\square$}
\newcommand{\Remark}{\note{Remark}}
\newcommand{\modin}{$\,$\\[-4mm] \indent}
\newcommand{\forevery}{\quad \forall}
\newcommand{\set}[1]{\{#1\}}
\newcommand{\setdef}[2]{\{\,#1:\,#2\,\}}
\newcommand{\setm}[2]{\{\,#1\mid #2\,\}}
\newcommand{\mt}{\mapsto}
\newcommand{\lra}{\longrightarrow}
\newcommand{\lla}{\longleftarrow}
\newcommand{\llra}{\longleftrightarrow}
\newcommand{\Lra}{\Longrightarrow}
\newcommand{\Lla}{\Longleftarrow}
\newcommand{\Llra}{\Longleftrightarrow}
\newcommand{\warrow}{\rightharpoonup}
\newcommand{
\paran}[1]{\left (#1 \right )}
\newcommand{\sqbr}[1]{\left [#1 \right ]}
\newcommand{\curlybr}[1]{\left \{#1 \right \}}
\newcommand{\abs}[1]{\left |#1\right |}
\newcommand{\norm}[1]{\left \|#1\right \|}
\newcommand{
\paranb}[1]{\big (#1 \big )}
\newcommand{\lsqbrb}[1]{\big [#1 \big ]}
\newcommand{\lcurlybrb}[1]{\big \{#1 \big \}}
\newcommand{\absb}[1]{\big |#1\big |}
\newcommand{\normb}[1]{\big \|#1\big \|}
\newcommand{
\paranB}[1]{\Big (#1 \Big )}
\newcommand{\absB}[1]{\Big |#1\Big |}
\newcommand{\normB}[1]{\Big \|#1\Big \|}
\newcommand{\produal}[1]{\langle #1 \rangle}

\newcommand{\thkl}{\rule[-.5mm]{.3mm}{3mm}}
\newcommand{\thknorm}[1]{\thkl #1 \thkl\,}
\newcommand{\trinorm}[1]{|\!|\!| #1 |\!|\!|\,}
\newcommand{\bang}[1]{\langle #1 \rangle}
\def\angb<#1>{\langle #1 \rangle}
\newcommand{\vstrut}[1]{\rule{0mm}{#1}}
\newcommand{\rec}[1]{\frac{1}{#1}}
\newcommand{\opname}[1]{\mbox{\rm #1}\,}
\newcommand{\supp}{\opname{supp}}
\newcommand{\dist}{\opname{dist}}
\newcommand{\myfrac}[2]{{\displaystyle \frac{#1}{#2} }}
\newcommand{\myint}[2]{{\displaystyle \int_{#1}^{#2}}}
\newcommand{\mysum}[2]{{\displaystyle \sum_{#1}^{#2}}}
\newcommand {\dint}{{\displaystyle \myint\!\!\myint}}
\newcommand{\q}{\quad}
\newcommand{\qq}{\qquad}
\newcommand{\hsp}[1]{\hspace{#1mm}}
\newcommand{\vsp}[1]{\vspace{#1mm}}
\newcommand{\ity}{\infty}
\newcommand{\prt}{\partial}
\newcommand{\sms}{\setminus}
\newcommand{\ems}{\emptyset}
\newcommand{\ti}{\times}
\newcommand{\pr}{^\prime}
\newcommand{\ppr}{^{\prime\prime}}
\newcommand{\tl}{\tilde}
\newcommand{\sbs}{\subset}
\newcommand{\sbeq}{\subseteq}
\newcommand{\nind}{\noindent}
\newcommand{\ind}{\indent}
\newcommand{\ovl}{\overline}
\newcommand{\unl}{\underline}
\newcommand{\nin}{\not\in}
\newcommand{\pfrac}[2]{\genfrac{(}{)}{}{}{#1}{#2}}

\def\ga{\alpha}     \def\gb{\beta}       \def\gg{\gamma}
\def\gc{\chi}       \def\gd{\delta}      \def\ge{\epsilon}
\def\gth{\theta}                         \def\vge{\varepsilon}
\def\gf{\phi}       \def\vgf{\varphi}    \def\gh{\eta}
\def\gi{\iota}      \def\gk{\kappa}      \def\gl{\lambda}
\def\gm{\mu}        \def\gn{\nu}         \def\gp{\pi}
\def\vgp{\varpi}    \def\gr{\rho}        \def\vgr{\varrho}
\def\gs{\sigma}     \def\vgs{\varsigma}  \def\gt{\tau}
\def\gu{\upsilon}   \def\gv{\vartheta}   \def\gw{\omega}
\def\gx{\xi}        \def\gy{\psi}        \def\gz{\zeta}
\def\Gg{\Gamma}     \def\Gd{\Delta}      \def\Gf{\Phi}
\def\Gth{\Theta}
\def\Gl{\Lambda}    \def\Gs{\Sigma}      \def\Gp{\Pi}
\def\Gw{\Omega}     \def\Gx{\Xi}         \def\Gy{\Psi}

\def\CS{{\mathcal S}}   \def\CM{{\mathcal M}}   \def\CN{{\mathcal N}}
\def\CR{{\mathcal R}}   \def\CO{{\mathcal O}}   \def\CP{{\mathcal P}}
\def\CA{{\mathcal A}}   \def\CB{{\mathcal B}}   \def\CC{{\mathcal C}}
\def\CD{{\mathcal D}}   \def\CE{{\mathcal E}}   \def\CF{{\mathcal F}}
\def\CG{{\mathcal G}}   \def\CH{{\mathcal H}}   \def\CI{{\mathcal I}}
\def\CJ{{\mathcal J}}   \def\CK{{\mathcal K}}   \def\CL{{\mathcal L}}
\def\CT{{\mathcal T}}   \def\CU{{\mathcal U}}   \def\CV{{\mathcal V}}
\def\CZ{{\mathcal Z}}   \def\CX{{\mathcal X}}   \def\CY{{\mathcal Y}}
\def\CW{{\mathcal W}} \def\CQ{{\mathcal Q}}
\def\BBA {\mathbb A}   \def\BBb {\mathbb B}    \def\BBC {\mathbb C}
\def\BBD {\mathbb D}   \def\BBE {\mathbb E}    \def\BBF {\mathbb F}
\def\BBG {\mathbb G}   \def\BBH {\mathbb H}    \def\BBI {\mathbb I}
\def\BBJ {\mathbb J}   \def\BBK {\mathbb K}    \def\BBL {\mathbb L}
\def\BBM {\mathbb M}   \def\BBN {\mathbb N}    \def\BBO {\mathbb O}
\def\BBP {\mathbb P}   \def\BBR {\mathbb R}    \def\BBS {\mathbb S}
\def\BBT {\mathbb T}   \def\BBU {\mathbb U}    \def\BBV {\mathbb V}
\def\BBW {\mathbb W}   \def\BBX {\mathbb X}    \def\BBY {\mathbb Y}
\def\BBZ {\mathbb Z}

\def\GTA {\mathfrak A}   \def\GTB {\mathfrak B}    \def\GTC {\mathfrak C}
\def\GTD {\mathfrak D}   \def\GTE {\mathfrak E}    \def\GTF {\mathfrak F}
\def\GTG {\mathfrak G}   \def\GTH {\mathfrak H}    \def\GTI {\mathfrak I}
\def\GTJ {\mathfrak J}   \def\GTK {\mathfrak K}    \def\GTL {\mathfrak L}
\def\GTM {\mathfrak M}   \def\GTN {\mathfrak N}    \def\GTO {\mathfrak O}
\def\GTP {\mathfrak P}   \def\GTR {\mathfrak R}    \def\GTS {\mathfrak S}
\def\GTT {\mathfrak T}   \def\GTU {\mathfrak U}    \def\GTV {\mathfrak V}
\def\GTW {\mathfrak W}   \def\GTX {\mathfrak X}    \def\GTY {\mathfrak Y}
\def\GTZ {\mathfrak Z}   \def\GTQ {\mathfrak Q}

\font\Sym= msam10 
\def\SYM#1{\hbox{\Sym #1}}
\newcommand{\bdw}{\prt\Gw\xspace}
\date{}
\maketitle\medskip

\noindent{\small {\bf Abstract} We study local and global properties of positive solutions of $-\Gd u=u^p+M\abs{\nabla u}^q$ in a domain $\Gw$ of $\BBR^N$, in the range $\min\{p,q\}>1$ and $M\in\BBR$. We prove a priori estimates and existence or non-existence of ground states for the same equation.
}\smallskip

\tableofcontents
\vspace{1mm}
\hspace{.05in}
\medskip
\mysection{Introduction}
This article is concerned with local and global properties of positive solutions of the following type of equations
\bel{I-0'}
-\Gd u=M'|u|^{p-1}u+M\abs{\nabla u}^q,
\ee
in $\Gw\setminus\{0\}$ where $\Gw$ is an open subset of $\BBR^N$ containing $0$, $p$ and $q$ are exponents larger than $1$ and $M,M'$ are real parameters. If $M'\leq 0$ the equation satisfies a comparison principle and a big part of the study can be carried via radial local supersolutions. This no longer the case when $M'>0$ which will be assumed in all the article, and by homothety $(\ref{I-0'})$ becomes
\bel{I-0}
-\Gd u=|u|^{p-1}u+M\abs{\nabla u}^q.
\ee
If $M=0$ $(\ref{I-0})$ is called Lane-Emden equation
\bel{I-1}
-\Gd u=|u|^{p-1}u.
\ee
It turns out that it plays an important role in modelling meteorological or astrophysical phenomena \cite{Emd}, \cite{Chandr}, this is the reason for which the first study, in the radial case, goes back to the end of nineteenth century and the beginning of the twentieth. A fairly complete presentation can be found in \cite{Fowl}. If $N\geq 3$, This equations exhibits two main critical exponents $p=\frac{N}{N-2}$ and $p=\frac{N+2}{N-2}$ which play a key role in the description of the set of positive solutions which can be summarized by the following overview:
 \smallskip

\nind 1- If $1<p\leq\frac{N}{N-2}$, there exists no positive solution if $\Gw$ is the complement of a compact set. Even in that case solution can be replaced by supersolution. This is easy to prove by studying the inequality satisfied by the spherical average of a solution of the equation.\smallskip

\nind 2- If $1<p<\frac{N+2}{N-2}$, there exists no {\it ground state}, i.e. positive solution in $\BBR^N$. Furthermore any positive solution $u$ in a ball $B_R=B_R(a)$ satisfies
\bel{I-2}
u(x)\leq c(R-|x-a|)^{-\frac{2}{p-1}},
\ee
 where $c=c(N,p)>0$, see \cite{GiSp}. \smallskip

\nind 3- If $p=\frac{N+2}{N-2}$ all the positive solutions in $\BBR^N$ are radial with respect to some point $a$ and endow the following form
\bel{I-3}
u(x):=u_\gl(x)=\frac{(N(N-2)\gl)^{\frac{N-2}{4}}}{\left(\gl+|x-a|^{2}\right)^{\frac{N-2}{2}}}.
\ee
All the positive solutions in $\BBR^N\setminus\{0\}$ are radial, see \cite{CaGiSp}. \smallskip

\nind 4- If $p>\frac{N+2}{N-2}$ there exist infinitely many positive ground states radial with respect to some points. They are obtained from one
say $v$, radial for example with respect to $0$ by the scaling transformation $T_k$ where $k>0$ with
\bel{I-4}
T_k[v](x)=k^{\frac2{p-1}}v(kx).
\ee

Indeed, the first significant non-radial results deals with  the case $1<p\leq\frac{N}{N-2}$. They are based upon the Brezis-Lions lemma \cite{BL} which yields an estimate of solutions in the Lorentz space $L^{\frac{N}{N-2},\infty}$, implying in turn the local integrability of
$u^q$. Then a bootstrapping method as in \cite{Lio} leads easily to some a priori estimate. Note that this subcritical case can be interpreted using the famous Serrin's results on quasilinear equations \cite{Se1}. The first breakthrough in the study of Lane-Emden equation came
in the treatment of the case $1<p<\frac{N+2}{N-2}$; it is due to Gidas and Spruck \cite{GiSp}. Their analysis is based upon differentiating the equation and then obtaining sharp enough local integral estimates on the term $u^{q-1}$ making possible the utilization of Harnack inequality as in \cite{Se1}. The treatment of the critical case $p=\frac{N+2}{N-2}$, due to Caffarelli, Gidas and Spruck \cite {CaGiSp}, was made possible thanks to a completely new approach based upon a combination of moving plane analysis and geometric measure theory. As for the supercritical case, not much is known and the existence of radial ground states is a consequence of Pohozaev's identity \cite{Po}, using a shooting method. \smallskip

The study of $(\ref{I-0})$ when $M\neq 0$ presents some similarities with the one of Lane-Emden equation in the cases 1 and 2, except that the proof are much more involved. Actually the approach we develop in this article is much indebted to our recent paper \cite{BiGaVe1} where we study local and global aspects of positive solutions of
\bel{I-01}
-\Gd u=u^p\abs{\nabla u}^q,
\ee
where $p\geq 0$, $0\leq q<2$, mostly in the superlinear case $p+q-1>0$. Therein we prove the existence of a critical line of exponents
\bel{I-02}
(\frak L):=\{(p,q)\in \BBR_+\ti [0,2): (N-2)p+(N-1)q=N\}.
\ee
The subcritical range corresponds to the fact that $(p,q)$ is below $(\frak L)$. In this region Serrin's celebrated results \cite{Se1} can be applied and we prove \cite[Theorem A]{BiGaVe1} that positive solutions of $(\ref{I-01})$ in the punctured ball $B_2\setminus\{0\}$ satisfy, for some constant $c>0$ depending on the solution,
\bel{I-03}
u(x)+\abs x\abs{\nabla u(x)}\leq c\abs x^{2-N}\quad\text{for all }x\in B_{1}\setminus\{0\}.
\ee

When $(p,q)$ is above $(\frak L)$, i.e. in the supercritical range, we introduced two methods for obtaining a priori estimate of solutions: The {\it pointwise Bernstein method} and the {\it integral Bernstein method}. The first one is based upon the change of unknown $u=v^{-\gb}$, and then to show that $\abs{\nabla v}$ satisfies an inequality of Keller-Osserman type. When $(p,q)$ lies above $(\frak L)$ and verifies \smallskip

\nind (i) either $1\leq p<\frac{N+3}{N-1}$ and $p+q-1<\frac{4}{N-1}$,\smallskip

\nind (ii) or $0\leq p<1$ and $p+q-1<\frac{(p+1)^2}{p(N-1)}$,\smallskip

\nind we prove that any positive solution of $(\ref{I-01})$ in a domain $\Gw\subset\BBR^N$ satisfies

\bel{I-04}
\abs {\nabla u^a(x)}\leq c^*\left(\dist (x,\prt\Gw)\right)^{-1-a\frac{2-q}{p+q-1}}\quad\text{for all }x\in \Gw,
\ee
for some positive $c^*$ and $a$ depending on $N$, $p$ and $q$ \cite[Theorem B]{BiGaVe1}. As a consequence we prove that any positive solution of $(\ref{I-01})$ in $\BBR^N$ is constant. With the second method we combine the change of unknown $u=v^{-\gb}$ with integration and cut-off functions. We show the existence of a quadratic polynomial $G$ in two variables such that  for any $(p,q)\in \BBR_+\ti [0,2)$ satisfying $G(p,q)<0$ any positive solution of $(\ref{I-01})$ in $\BBR^N$ is constant \cite[Theorem C]{BiGaVe1}. The polynomial $G$ is not simple but it is worth noting that if $0\leq p<\frac{N+2}{N-2}$, there holds $G(p,0)<0$, which recovers Gidas and Spruck result \cite{GiSp}.

\medskip

For equation $(\ref{I-0})$ we first observe that the equation is invariant under the scaling transformation $(\ref{I-4})$ for any $k>0$ if and only if
$q$ is {\it critical with respect to }$p$, i.e.
$$q=\frac{2p}{p+1}.$$
In general the transformation $T_k$ exchanges $(\ref{I-0})$ with
 \bel{I-5}
-\Gd v=v^p+Mk^{\frac{2p-q(p+1)}{p-1}}|\nabla v|^q,
\ee
hence if $q<\frac{2p}{p+1}$, the limit equation when $k\to 0$ is $(\ref{I-1})$. We say that the exponent $p$ is dominant. We can also consider the transformation
 \bel{I-6}
S_k[v](x)=k^{\frac{2-q}{q-1}}v(kx),
\ee
when $q\neq 2$, which is the same as $T_k$ if $q=\frac{2p}{p+1}$, and more generally transforms $(\ref{I-0})$ into
 \bel{I-7}
-\Gd v=k^{\frac{q-p(2-q)}{q-1}}v^p+M|\nabla v|^q.
\ee
Hence if $q>\frac{2p}{p+1}$, the limit equation when $k\to 0$ is the Riccati equation
 \bel{I-8}
-\Gd v=M|\nabla v|^q.
\ee
{\it It is also important to notice that the value of the coefficient $M$ (and not only its sign) plays a fundamental role, only if $q=\frac{2p}{p+1}$}. If $q\neq\frac{2p}{p+1}$ the transformation
\bel{I-9}
u(x)=av(y)\quad\text{with }\, a=\abs{M}^{-\frac{2}{(p+1)q-2p}}\text{ and }\, y=a^{\frac{p-1}{2}}x
\ee
allows to transform $(\ref{I-0})$ into
\bel{I-10}
-\Gd v=|v|^{p-1}v\pm|\nabla v|^q.
\ee

The equation $(\ref{I-0})$ has been essentially studied in the {\it radial case} when $M<0$ in connection with the parabolic equation
\bel{I-11}
\prt_tu-\Gd u+M|\nabla u|^q=|u|^{p-1}u,
\ee
see \cite{ChWe}, \cite{Fi}, \cite{FiQu}, \cite{SeZo}, \cite{So}, \cite{Vo}, \cite{Vo2}. The studies mainly deal with the case  $q\neq\frac{2p}{p+1}$, although not complete when $q>\frac{2p}{p+1}$. When $q=\frac{2p}{p+1}$ the existence of a ground state is proved in dimension 1. Some partial results that we will improve, already exist in higher dimension. The case $M>0$ attracted less attention.\smallskip

In the {\it nonradial case}, any nonnegative nontrivial solution is positive since $p,q>1$. We first observe, using a standard averaging method applied to positive supersolutions of ($\ref{I-1})$, that if $M\geq 0$, $1<p\leq \frac{N}{N-2}$ when $N\geq 3$, any $p>1$ if $N=1,2$, then for any $q>0$ there exists no positive solution in an exterior domain. When $0<q<\frac{2p}{p+1}$ the equation endows some character of the pure Emden-Fowler equation $(\ref{I-1})$ by the transformation $T_k$. In \cite{PoQuSo} it is proved that if $0<q<\frac{2p}{p+1}$, $1<p<\frac{N+2}{N-2}$ and $M\in\BBR$, any positive solution of $(\ref{I-1})$ in an open domain satisfies
\bel{I-12}
u(x)+\abs{\nabla u(x)}^{\frac{2}{p+1}}\leq c_{N,p,q,M}\left(1+\left(\dist (x,\prt\Gw)\right)^{-\frac{2}{p-1}}\right)\qquad\text{for all }x\in\Gw.
\ee
Note that this does not imply the non-existence of ground state. In \cite{AGMQ1} Alarc\'on, Garc\'ia-Meli\'an and Quass study the equation
\bel{I-12-1}
-\Gd u=|\nabla u|^q+f(u),
\ee
in an exterior domain of $\BBR^N$ emphasizing the fact that positive solutions are super harmonic functions. They prove that if $1<q\leq \frac{N}{N-1}$ and if $f$ is positive on $(0,\infty)$ and satisfies
\bel{I-12-2}\displaystyle
\limsup_{s\to 0}s^{-p}f(s)>0,
\ee
for some $p> \frac{N}{N-2}$, then $(\ref{I-12-1})$ admits no positive supersolution. The same authors also study in \cite{AGMQ2} existence and non-existence of positive solutions of $(\ref{I-12-1})$ in a bounded domain with Dirichlet condition.
\smallskip

The techniques we developed in this paper are based upon a delicate extension of the ones already introduced in \cite{BiGaVe1}. Our first nonradial result dealing with the case $q>\frac{2p}{p+1}$ is the following:
\medskip

\nind{\bf Theorem A} {\it Let $N\geq 1$, $p>1$ and $q>\frac{2p}{p+1}$. Then for any $M>0$, any solution of $(\ref{I-0})$ in a domain $\Gw\subset\BBR^N$ satisfies
\bel{I-13}
\abs{\nabla u(x)}\leq c_{N,p,q}\left(M^{-\frac{p+1}{(p+1)q-2p}}+\left(M\dist (x,\prt\Gw)\right)^{-\frac{1}{q-1}}\right)\quad\text{for all }x\in\Gw.
\ee
As a consequence, any ground state has at most a linear growth at infinity:
 \bel{I-13+}
\abs{\nabla u(x)}\leq c_{N,p,q}M^{-\frac{p+1}{(p+1)q-2p}}\qquad\text{for all }x\in\BBR^N.
\ee
}

Our proof relies on a direct Bernstein method combined with Keller-Osserman's estimate applied to $\abs{\nabla u}^2$. It is important to notice that the result holds for any $p>1$, showing that, in some sense, the presence of the gradient term has a regularizing effect. In the case $q<\frac{2p}{p+1}$ we prove a non-existence result\medskip

\nind{\bf Theorem A'} {\it Let $N\geq 1$, $p>1$, $1<q<\frac{2p}{p+1}$ and $M>0$. Then there exists a constant $c_{N,p,q}>0$ such that there is no positive solution of $(\ref{I-0})$ in $\BBR^N$ satisfying
\bel{I-13-1}
u(x)\leq c_{N,p,q}M^{\frac{2}{2p-(p+1)q}}\qquad\text{for all }x\in\BBR^N.
\ee}
\medskip

When $q$ is critical with respect to $p$ the situation is more delicate since the value of $M$ plays a fundamental role. Our first statement is a particular case of a more general result in \cite{AGMQ1}, but with a simpler proof which allows us to introduce techniques that we use later on. \medskip

\nind{\bf Theorem B} {\it Let $N\geq 2$, $p>1$ if $N=2$ or $1<p\leq \frac{N}{N-2}$ if $N=3$, $q=\frac{2p}{p+1}$ and $M>-\gm^*$ where
\bel{I-14}
\gm^*:=\gm^*(N)=(p+1)\left(\myfrac{N-(N-2)p}{2p}\right)^{\frac p{p+1}}.
\ee
Then there exists no nontrivial nonnegative supersolution of $(\ref{I-0})$ in an exterior domain.
}\medskip

In this range of values of $p$ this result is optimal since for $M\leq -\gm^*$ there exists positive singular solutions. The constant $\gm^*$ will play an important role in the description developed in \cite{BiGaVe2} of radial solutions of $(\ref{I-0})$. Using a variant of the method used in the proof of Theorem B we obtain  results of existence and nonexistence of large solutions.

\medskip

\nind{\bf Theorem B'} {\it Let $N\geq 1$, $p>1$ and $q=\frac{2p}{p+1}$.\smallskip

\nind1- If $\Gw$ is a domain with a compact boundary satisfying the Wiener criterion and
$ M\geq -\gm^*(2)$ there exists
no positive supersolution of $(\ref{I-0})$ in $\Gw$ satisfying
\bel{I-14'}\displaystyle
\lim_{\dist(x,\prt\Gw)\to 0}u(x)=\infty.
\ee

\nind 2- If $G$ is a bounded convex domain, $\Gw=\overline G^c$ and $ M< -\gm^*(1)$ there exists a positive solution of $(\ref{I-0})$ in $\Gw$ satisfying $(\ref{I-14'})$.
}\medskip

We show in \cite{BiGaVe2} that the inequality $ M< -\gm^*(1)$  is the necessary and sufficient condition for the existence of a radial large solution in the exterior of a ball. \smallskip

Concerning ground states, we prove their nonexistence for any $p>1$ provided $M>0$ is large enough: indeed \medskip

\nind{\bf Theorem C} {\it Let $\Gw\subset\BBR^N$, $N\geq 1$, be a domain, $p>1$, $q=\frac{2p}{p+1}$. For any
\bel{I-14b}M>M_\dag:=\left(\frac{p-1}{p+1}\right)^{\frac{p-1}{p+1}}\left(\frac{N(p+1)^2}{4p}\right)^{\frac{p}{p+1}},
\ee
and any $\gn>0$ such that $(1-\gn)M>M_\dag$, there exists a positive constant $c_{N,p,\gn}$ such that any solution $u$ in $\Gw$ satisfies
 \bel{I-12b}
\abs{\nabla u(x)}\leq c_{N,p,\gn}\left((1-\gn)M-M_\dag\right)^{-\frac{p+1}{p-1}}\left(\dist (x,\prt\Gw)\right)^{-\frac{p+1}{p-1}}\quad\text{for all }x\in\Gw.
\ee
Consequently there exists no nontrivial solution of $(\ref{I-0})$ in $\BBR^N$.
}\medskip

The next result, based upon an elaborate Bernstein method, complements Theorem C under a less restrictive assumption on $M$  but a more restrictive assumption on $p$.\medskip

\nind{\bf Theorem D} {\it Let  $1<p<\frac{N+3}{N-1}$, $N\geq 2$, $1<q<\frac{N+2}{N}$ and  $\Gw\subset\BBR^N$ be a domain.  Then  there exist $a>0$ and $c_{N,p,q}>0$ such that  for any $M>0$, any positive solution $u$ in $\Gw$ satisfies
 \bel{I-12c}
\abs{\nabla u^{a}(x)}\leq c_{N,p,q}\left(\dist(x,\prt\Gw)\right)^{-\frac{2a}{p-1}-1}\quad\text{for all }x\in\Gw.
\ee
Hence there exists no nontrivial nonnegative solution of $(\ref{I-0})$ in $\BBR^N$.
}\medskip

It is remarkable that the constants $a$ and $c_{N,p,q}$ do not depend on $M>0$, a fact which is clear when $q\neq\frac{2p}{p+1}$ by using the transformation $T_k$, but much more delicate to highlight when $q=\frac{2p}{p+1}$ since $(\ref{I-0})$ is invariant. When $\abs M$ is small, we use an integral method to obtain the following result which contains, as a particular case, the estimates in \cite {GiSp} and \cite {BiGaVe2}. The key point of this method is to prove that the solutions in a punctured domain satisfy a local Harnack inequality.  \medskip

\nind{\bf Theorem E} {\it Let $N\geq 3$,  $1<p<\frac{N+2}{N-2}$, $q=\frac{2p}{p+1}$.  Then there exists $\ge_0>0$ depending on $N$ and $p$ such that for any $M$ satisfying $\abs M\leq\ge_0$, any positive solution $u$ in $B_R\setminus\{0\}$ satisfies
 \bel{I-12x}
u(x)\leq c_{N,p}\abs x^{-\frac{2}{p-1}}\quad\text{for all }x\in B_{\frac R2}\setminus\{0\}.
\ee
As a consequence there exists no positive solution of $(\ref{I-0})$ in $\BBR^N$, and any positive solution $u$ in a domain $\Gw$ satisfies
 \bel{I-12xx}
u(x)+\abs{\nabla u(x) }^{\frac{2}{p+1}}\leq c'_{N,p}\left(\dist(x,\prt\Gw)\right)^{-\frac{2}{p-1}}\quad\text{for all }x\in \Gw.
\ee
}\medskip

Note that under the assumptions of Theorem E, there exist ground states for $\abs M$ large enough when $1<p<\frac{N}{N-2}$, or any
$p>1$ if $N=1,2$. \medskip

If $u$ is a radial  solutions of $(\ref{I-0})$ in $\BBR^N$ it satisfies
 \bel{I-13x}
-u''-\myfrac{N-1}{r}u'=\abs u^{p-1}u+M\abs {u'}^{q},
\ee
on $(0,\infty)$. Using several type of Lyapounov type functions introduced  by Leighton \cite{Lei} and Anderson and Leighton \cite{AnLe}, we prove some results dealing with the case $M>0$ which complement the ones  of \cite{SeZo} relative to the case $M<0$. \medskip

\nind{\bf Theorem F} {\it 1- Let $p>1$ and $q>\frac{2p}{p+1}$. Then there exists no radial ground state $u$ satisfying $u(0)=1$ when $M>0$ is too large.
\smallskip

\nind 2- Let $1<p<\frac{N+2}{N-2}$. If $1<q\leq p$ there exists no radial ground state for any $M>0$. If $q>p$
there exists no radial ground state for $M>0$ small enough.\smallskip

\nind 3- Let $N\geq 3$, $p>\frac{N+2}{N-2}$ and $q\geq \frac{2p}{p+1}$. Then there exist radial ground states for $M>0$ small enough.
}\medskip

We end the article in proving the existence of non-radial positive singular solutions of $(\ref{I-0})$ in $\BBR^N\setminus\{0\}$ in the case $q=\frac{2p}{p+1}$ obtained by bifurcation from radial explicit positive singular solutions.
Our result shows that the situation is very contrasted according $M>0$ where a bifurcation from $(M,X_{M})$ occurs only if $p\geq\frac{N+1}{N-3}$ and $M\geq 0$  and $M<0$ where there exists a countable set of bifurcations from $(M_k,X_{M_k})$, $k\geq 1$, when $1<p<\frac{N+1}{N-3}$.
\medskip

In a subsequent article \cite{BiGaVe2} we present a fairly complete description of the positive radial solutions of $(\ref{I-0})$ in $\BBR^N\setminus\{0\}$ in the scaling invariant case $q=\frac{2p}{p+1}$.\medskip

\nind{\bf Acknowledgements} This article has been prepared with the support of the collaboration
programs ECOS C14E08 and FONDECYT grant 1160540 for the three authors.

\mysection{The direct Bernstein method}
We begin with a simple property in the case $M\geq 0$ which is a consequence of the fact that the positive solutions of $(\ref{I-0})$ are superharmonic.
\bprop{zero} 1- There exists no positive solution of $(\ref{I-0})$ in $\BBR^N\setminus\overline B_R$, $R\geq 0$ if one of the two conditions is satisfied:\\
(i) $M\geq 0$, $q\geq 0$ and either $N=1,2$ and $p>1$ or $N\geq 3$ and $1<p\leq \frac{N}{N-2}$.\\
(ii) $M> 0$, $N\geq 3$, $p\geq 1$ and $1<q\leq \frac{N}{N-1}$.\smallskip

\nind 2- If $N\geq 3$, $q\geq1$, $p>\frac{N}{N-2}$ and $u(x)=u(r,\gs)$ is a positive solution of $(\ref{I-0})$ in $\BBR^N\setminus\overline B_R$, $R\geq 0$. Then  there exists $\gr\geq R$ such that
\bel{I-1-01}
\myfrac{1}{N\gw_N}\myint{S^{N-1}}{}u(r,\gs)dS:=\overline u(r)\leq c_0r^{-\frac{2}{p-1}}\qquad\text{for all }\;r>\gr ,
\ee
with $c_0:=\left(\frac{2N}{p-1}\right)^{\frac{1}{p-1}}$ and
\bel{I-1-01'}
\abs{\myfrac{1}{N\gw_N}\myint{S^{N-1}}{}u_r(r,\gs)dS}:=\abs{\overline u_r(r)}\leq (N-2)c_0r^{-\frac{p+1}{p-1}}\qquad\text{for all }\;r>\gr.
\ee
\nind 3- If $M>0$, $p\geq 0$, and $q>\frac{N}{N-1}$ there holds for
\bel{I-1-01''}
\abs{\overline u_r(r)}\leq \left(\myfrac{(q-1)(N-1)-1}{(q-1)M}\right)^{\frac{1}{q-1}}r^{-\frac{1}{q-1}}\qquad\text{for all }\;r>\gr,
\ee
and
\bel{I-1-01'''}
\overline u(r)\leq \Bigl(\myfrac{q-1}{2-q}\Bigr)\left(\myfrac{(q-1)(N-1)-1}{(q-1)M}\right)^{\frac{1}{q-1}}r^{\frac{q-2}{q-1}}\qquad\text{for all }\;r>\gr,
\ee
Furthermore, if $R=0$, inequalities $(\ref{I-1-01})$,  $(\ref{I-1-01'})$ and $(\ref{I-1-01''})$ hold with $\gr=0$.
\es
\Proof Assertion 1-(i) is not difficult to obtain by integrating the inequality satisfied by the spherical average of the solution and using Jensen's inequality. For the sake of completeness, we give a simple proof although the result is actually valid for much more general equations (see e.g. \cite{BiPo} and references therein). In this statement we denote by $(r,\gs)\in\BBR_+\ti S^{N-1}$ the spherical coordinates in $\BBR^N$, by $\gw_N$ the volume of the unit N-ball and thus $N\gw_N$ is the (N-1)-volume of the unit sphere $S^{N-1}$.
Writing $(\ref{I-0})$ in spherical coordinates and using Jensen formula, we get
\bel{I-1-02}-r^{1-N}\left(r^{N-1}\overline u_r\right)_r\geq \overline u^p+M\abs{\overline u_r}^q.
\ee
It implies that $r\mapsto w(r):= -r^{N-1}\overline u_r$ is increasing on $(R,\infty)$, thus it admits a limit $\ell\in(-\infty,\infty]$. If $\ell\leq 0$, then $\overline u_r(r)> 0$ on $(R,\infty)$. Hence $\overline u(r)\geq\overline u(\gr):=c>0$ for $r\geq \gr>R$. then
$$\left(r^{N-1}\overline u_r\right)_r\leq -c^pr^{N-1}\Longrightarrow \overline u_r(r)\leq \myfrac{\gr^{N-1}}{r^{N-1}}\overline u_r(\gr)
-\myfrac{c^p}{N}\left(r-\myfrac{\gr^{N}}{r^{N-1}}\right),
$$
which implies $ \overline u_r(r)\to-\infty$, thus $ \overline u(r)\to-\infty$ as $r\to-\infty$, contradiction. Therefore $\ell\in(0,\infty]$ and either
$\overline u_r(r)<0$ on $(R,\infty)$ or there exists $r_{\ell}>R$ such that $\overline u_r(r_\ell)=0$, $\overline u$ is increasing on $(R,r_\ell,)$ and decreasing on $(r_\ell,\infty)$. If $\overline u_r(r)<0$ on $(R,\infty)$, then we have for $r> 2R$
$$-r^{N-1}\overline u_r(r)\geq \myint{\frac {r}{2}}{r}t^{N-1}\overline u^p(t)dt\geq \myfrac{r^{N}\overline u^p(r)}{2N}
\Longrightarrow \left(\overline u^{1-p}\right)_r\geq \myfrac{(p-1)r}{2N}\Longrightarrow\overline u(r)\leq \left(\myfrac{2N}{(p-1)r^2}\right)^{\frac{1}{p-1}},
$$
which yields $(\ref{I-1-01})$.
If we are in the second case with $r_\ell>R$, we apply the same inequality with $r> 2r_\ell$ and again $(\ref{I-1-01})$ for $r> 2r_\ell$.
Since  $\overline u$ is superharmonic, the function $v(s)=\overline u(r)$ with $s=r^{2-N}$ is concave on
$(0,R^{2-N})$ and it tends to $0$ when $s\to 0$. Thus
$$v_s(s)\leq \frac{v}{s}\Longrightarrow \abs{\overline u_r(r)}\leq (N-2)\myfrac{\overline u(r)}{r}\leq (N-2)c_0r^{-\frac{p+1}{p-1}}.$$
This implies $(\ref{I-1-01})$ and $(\ref{I-1-01'})$. Note that the case $r_\ell>R$ cannot happen if $R=0$, so in any case, if $R=0$ then $\gr=0$.\\
If $M>0$, we have with $w(r)=-r^{N-1}\overline u_r$
$$w_r\geq Mr^{(1-q)(N-1)}\abs{w}^q.
$$
We have seen that $w(r)>0$ at infinity with limit $\ell\in (0,\infty]$, hence, on the maximal interval containing $\infty$ where $w>0$, we have $(w^{1-q})_r\leq (1-q)Mr^{(N-1)(1-q)}$. We have for $r>s>R$
$$w^{1-q}(r)-w^{1-q}(s)\leq M\ln\left(\frac{r}{s}\right),
$$
if $q=\frac{N}{N-1}$ and
$$w^{1-q}(r)-w^{1-q}(s)\leq \frac{M(q-1)}{(q-1)(N-1)-1}\left(r^{1-(q-1)(N-1)}-s^{1-(q-1)(N-1)}\right)
$$
if $q<\frac{N}{N-1}$, and both expressions which tend to $-\infty$ when $r\to\infty$, a contradiction. This proves 1-(ii).
If $q>\frac{N}{N-1}$, the above expression yields, when $r\to\infty$,
$$\ell^{1-q}-w^{1-q}(s)\leq -\myfrac{(q-1)M}{(q-1)(N-1)-1}s^{1-(q-1)(N-1)}.
$$
This implies
$$w(s)\leq \left(\myfrac{(q-1)(N-1)-1}{(q-1)M}\right)^{\frac{1}{q-1}}s^{N-1-\frac{1}{q-1}},
$$
and $(\ref{I-1-01''})$. $\phantom{---}$
$\phantom{---}$\qeda\medskip

\nind \Remark The previous is a particular case of a much more general one dealing with quasilinear operators proved in \cite[Theorem 3.1] {BiPo}.

\subsection{Proof of Theorems A, A' and C}
The function $u$ is at least $C^{3+\ga}$ for some $\ga\in (0,1)$ since $p,q>1$. Hence $z=\abs{\nabla u}^2$ is $C^{2+\ga}$. Since there holds by Bochner's identity and Schwarz's inequality
\bel{I-1-1}
-\myfrac{1}{2}\Gd z+\myfrac {1}{N}(\Gd u)^2+\langle\nabla\Gd u,\nabla u\rangle\leq 0,
\ee
we obtain from $(\ref{I-0})$,
$$
-\myfrac{1}{2}\Gd z+\myfrac {|u|^{2p}}{N}+\myfrac {2M}{N}|u|^{p-1}uz^{\frac q2}+\myfrac {M^2}{N}z^q-p|u|^{p-1}z
-\myfrac {Mq}{2}z^{\frac q2-1}\langle\nabla z,\nabla u\rangle\leq 0.
$$
Since for $\gd>0$,
$$z^{\frac q2-1}\left|\langle\nabla z,\nabla u\rangle\right|\leq \abs{z^{-\frac12}\nabla z}z^{\frac{q-1}2}\abs{\nabla u}=
\abs{z^{-\frac12}\nabla z}z^{\frac{q}2}\leq \gd z^q+\myfrac1{4\gd}\myfrac{\abs{\nabla z}^2}{z},
$$
we obtain for any $\gn\in (0,1)$, provided $\gd$ is small enough,
\bel{I-1-2}
-\myfrac{1}{2}\Gd z+\myfrac {|u|^{2p}}{N}+\myfrac {2M}{N}|u|^{p-1}uz^{\frac q2}+\myfrac {M^2(1-\gn)^2}{N}z^q-p|u|^{p-1}z
\leq c_1\myfrac{\abs{\nabla z}^2}{z},
\ee
where $c_1=c_1(M,N,\gn)>0$. \smallskip

\subsubsection {Proof of Theorem A}
We recall the following technical result proved in \cite[Lemma 2.2]{BiGaVe1} which will be used several times in the course of this article.

\blemma{L2.2} Let $S>1$, $R>0$ and $v$ be  continuous and nonnegative in $\overline B_R$ and $C^1$ on the set
$\CU_+=\{x\in B_R:v(x)>0\}$. If $v$ satisfies, for some real number $a$,
\bel{5-2.2}
-\Gd v+v^S\leq a\myfrac{|\nabla v|^2}{v}
\ee
on each connected component of $\CU_+$, then
\bel{5-2.2+}
v(0)\leq c_{N,S,a}R^{-\frac{2}{S-1}}.
\ee
\es

\nind {\it Abridged proof}. Assuming $a>0$, we set $W=v^\ga$ for $0<\ga\leq \frac{1}{a+1}$,
this transforms $(\ref{5-2.2})$ into
\bel{5-2.2++}
-\Gd W+\myfrac{1}{\ga}W^{\ga(S-1)+1}\leq 0,
\ee
and then we apply Keller-Osserman inequality.\qeda

\medskip

\nind{\it Proof of Theorem A.} Suppose $\frac{2p}{p+1}<q$. We set $r=\frac{2p}{p-1}$, $r'=\frac{r}{r-1}$, then, for any $\ge>0$
$$p|u|^{p-1}z\leq \myfrac{\ge^r|u|^{(p-1)r}}{r}+\myfrac{z^{r'}}{\ge^{r'}r'}=(p-1)\myfrac{\ge^r|u|^{2p}}{2}
+(p+1)\myfrac{z^{\frac{2p}{p+1}}}{2\ge^{r'}}.
$$
We fix $\eta\in (0,1)$ and $\ge$ so that $\ge^r=\frac{2(1-\eta)}{N(p-1)}$ and get
$$p|u|^{p-1}z\leq(1-\eta)\frac{|u|^{2p}}{N}+c_2z^{\frac{2p}{p+1}},
$$
where $c_2=\frac{p+1}{2}\left(\frac{N(p-1)}{2(1-\eta)}\right)^{\frac{p+1}{p-1}}$. We perform the change of scale $(\ref{I-4})$
in order to reduce  $(\ref{I-0})$ to the case $M=1$ by setting $u(x)=\ga^{\frac{2}{p-1}}v(\ga x)$ with $\ga=M^{-\frac{p-1}{(p+1)q-2p}}$. Then the equation for $z=\abs{\nabla v}^2$ is considered in $\Gw_\ga=\ga \Gw$. Choosing now $\eta=\frac12$ we obtain
$$c_2z^{\frac{2p}{p+1}}\leq \myfrac{1}{4N}z^q+c_3,
$$
where $c_3=c_3(N,p,q)>0$, hence
$$-\myfrac{1}{2}\Gd z+\myfrac{v^{2p}}{2N}+\myfrac{1}{4N}z^q\leq c_3+c_1\myfrac{\abs{\nabla z}^2}{z}.
$$
Put $\tilde z=\left(z-\left(4Nc_3\right)^{\frac 1q}\right)_+$, then
$$-\myfrac{1}{2}\Gd \tilde z+\myfrac{1}{4N}\tilde z^q\leq c_1\myfrac{\abs{\nabla \tilde z}^2}{\tilde z},
$$
hence, from \rlemma{L2.2}, we derive
$$\tilde z(y)\leq c_4\left(\dist (y,\prt\Gw_\ga)\right)^{\frac 2{q-1}}
$$
where $c_4=c_4(N,q,c_1)>0$ which implies
\bel{I-1-2-1}
\abs{\nabla v(y)}\leq c'_4\left(1+\left(\dist (y,\prt\Gw_\ga)\right)^{-\frac{1}{q-1}}\right)\qquad\forall\, y\in \Gw_\ga.
\ee
Then $(\ref{I-13})$ and $(\ref{I-13+})$ follow.\smallskip

Assume now that there exists a ground state $u$. Fix $y\in\BBR^N$ and consider $\{y_n\}\subset \BBR^N$ such that $|y_n|=2n>|y|$. We apply $(\ref{I-1-2-1})$ with $\Gw_\ga=B_n(y_n)$. Then
$$\abs{\nabla v(y)}\leq c'_4\left(1+\abs{2n-|y|}^{-\frac{1}{q-1}}\right),
$$
and letting $n\to\infty$ we infer
\bel{I-1-2-2}
\abs{\nabla v(y)}\leq c'_4\qquad\forall\, y\in \BBR^N.
\ee
Hence, by the definition of $v$ and $y$ we see that
$$|\nabla u(x)|\le c'_4M^{-\frac{p+1}{(p+1)q-2p}}\qquad\forall\, x\in \BBR^N$$
which is exactly \eqref{I-13+}.
\qeda

\medskip

\subsubsection{Proof of Theorem A'}
Suppose $1<q<\frac{2p}{p+1}$. By scaling we reduce to the case $M=1$ and we replace $u$ by $v$ defined by $(\ref{I-4})$ as in the proof of Theorem A with $\ga=M^{\frac{p-1}{2p-(p+1)q}}$. From $(\ref{I-1-2})$ with $\gn=\frac 14$ the function $z=\abs{\nabla v}^2$ satisfies

\bel{I-1-2-3}
-\myfrac{1}{2}\Gd z+\myfrac {v^{2p}}{N}+\myfrac {1}{2N}z^q-pv^{p-1}z
\leq c_1\myfrac{\abs{\nabla z}^2}{z}.
\ee
By H\"older's inequality,
$$pv^{p-1}z\leq \myfrac {1}{4N}z^q+p(4Np)^{q'-1}v^{(p-1)q'}. 
$$
Since $(p-1)q'=2p+\frac{2p-(p+1)q}{q-1}$ we derive
$$-\myfrac{1}{2}\Gd z+\myfrac {v^{2p}}{N}\left(1-4^{q'-1}p^{q'}N^{q'}v^{\frac{2p-(p+1)q}{q-1}}\right)+\myfrac {1}{4N}z^q\leq c_1\myfrac{\abs{\nabla z}^2}{z}.
$$
If $\max v\leq c_{N,p,q}:=(4^{q'-1}p^{q'}N^{q'})^{-\frac{q-1}{2p-(p+1)q}}$, we obtain
$$-\myfrac{1}{2}\Gd z+\myfrac {1}{4N}z^q\leq c_1\myfrac{\abs{\nabla z}^2}{z},
$$
which implies that $z=0$ by \rlemma{L2.2}, hence $v$ is constant and thus $v=0$ from the equation. $\phantom{---}$
$\phantom{---}$\qeda\medskip

\nind\Remark If $u$ is a positive ground state of $(\ref{I-0})$ radial with respect to $0$, it satisfies $u_r(0)=0$ and it is a decreasing function of $r$. The previous theorem asserts that it must satisfy
\bel{I-1-2-3'}
u(0)>c_{N,p,q} M^{\frac{2}{2p-(p+1)q}}.
\ee

\medskip

\subsubsection{Proof of Theorem C}
Suppose $\frac{2p}{p+1}=q$. For $A>0$ we  consider the expression 
$$\BA {lll}\left(u^p+A\abs{\nabla u}^{q}\right)^2-Npu^{p-1}\abs{\nabla u}^2\\[2mm]
\phantom{------}=\left(u^p+A\abs{\nabla u}^{q}-\sqrt{Np\,}u^{\frac{p-1}{2}}\abs{\nabla u}\right)\left(u^p+A\abs{\nabla u}^{q}+\sqrt{Np\,}u^{\frac{p-1}{2}}\abs{\nabla u}\right).
\EA$$
Now the function $Z\mapsto \Phi_A(Z)=u^{p}+AZ^{q}-\sqrt{Np\,}\,u^{\frac{p-1}{2}}Z$ achieves its minimum at
$Z_0=\left(\frac{\sqrt {Np}}{qA }\right)^{\frac{p+1}{p-1}}u^{\frac{p+1}{2}}$ and
$$\BA {lll}\Gf_A(Z_0)=\left[1-\myfrac{p-1}{p+1}\left(\myfrac{N(p+1)^2}{4p}\right)^{\frac{p}{p-1}}A^{-\frac{p+1}{p-1}}\right] u^{p}.
\EA$$
Thus setting
\bel{I-1-3}
M_\dag=\left(\frac{p-1}{p+1}\right)^{\frac{p-1}{p+1}}\left(\frac{N(p+1)^2}{4p}\right)^{\frac{p}{p+1}},
\ee
we obtain that if $A\geq M_\dag$, then $\Gf_A(Z)\geq 0$ for all $Z$. Put $M_\gn=(1-\gn) M$ for $\gn\in (0,1)$ such that $M_\dag<M_\gn$, we derive from
$(\ref{I-1-2})$
\bel{I-1-4}
-\myfrac{1}{2}\Gd z+\myfrac {(u^{p}+M_\dag z^{\frac q2})^2}{N}-pu^{p-1}z+\myfrac {M_\gn^2-M_\dag^2}{N}z^q
\leq c_1\myfrac{\abs{\nabla z}^2}{z},
\ee
which yields
$$-\myfrac{1}{2}\Gd z+\myfrac {M_\gn^2-M_\dag^2}{N}z^q
\leq c_1\myfrac{\abs{\nabla z}^2}{z}.
$$
Using again \rlemma{L2.2} we obtain
\bel{I-1-4b}
\abs{\nabla u(x)}\leq c'_1\left((1-\gn)M-M_\dag\right)^{-\frac{1}{q-1}}\left(\dist(x,\prt\Gw)\right)^{-\frac{1}{q-1}},
\ee
which is equivalent to $(\ref{I-12b})$.
\qeda

\subsection{Proof of Theorems B and B'}

\subsubsection{ Proof of Theorem B} Since the result is known when $M\geq 0$ from \rprop{zero}, we can assume that $M=-m<0$ and $N=1,2$ or $N\geq 3$ with $p<\frac{N}{N-2}$, $u$ is a nonnegative supersolution of $(\ref{I-0})$ in $\overline B^c_R$ and we set $u=v^b$ with $b>1$. Then
\bel{I-1-5}
-\Gd v\geq (b-1)\frac{\abs{\nabla v}^2}{v}+\frac{1}{b}v^{1+b(p-1)}-mb^{q-1}v^{(b-1)(q-1)}\abs{\nabla v}^q.
\ee
Here again $q=\frac{2p}{p+1}$, setting $z=\abs{\nabla v}^2$ we obtain
$$-\Gd v\geq \frac{\Phi(z)}{bv}
$$
where
$$\Gf(z)=b(b-1)z-mb^{\frac{2p}{p+1}}v^{\frac{2+b(p-1)}{p+1}}z^{\frac {p}{p+1}}+v^{2+b(p-1)}.$$
Thus $\Gf$ achieves it minimum for
$$z_0=\left(\frac{mpb^{q-1}}{(b-1)(p+1)}\right)^{p+1}b^{p-1}v^{2+b(p-1)}
$$
and
\bel{I-1-5*}\Gf(z_0)=v^{2+b(p-1)}\left(1-\frac {p^p}{(p+1)^{p+1}}\left(\frac{b}{b-1}\right)^pm^{p+1}\right).
\ee
In order to ensure the optimal choice, when $N\geq 3$ we take $1+b(p-1)=\frac{N}{N-2}$, hence $b=\frac{2}{(N-2)(p-1)}$ which is larger than $1$ because
$p<\frac{N}{N-2}$. Finally
$$
\Gf(z_0)=v^{\frac{N}{N-2}+1}\left(1-\frac {1}{(p+1)^{p+1}}\left(\frac{2p}{N-p(N-2)}\right)^pm^{p+1}\right).
$$
Hence, if
\bel{I-1-5'}
m< (p+1)\left(\frac{N-p(N-2)}{2p}\right)^{\frac{p}{p+1}}=\gm^*(N),
\ee
we have for some $\gd>0$,
\bel{I-1-5-1+}
-\Gd v\geq \gd v^{\frac{N}{N-2}},
\ee
and by \rprop{zero} that is no positive solution in an exterior domain of $\BBR^N$. \smallskip

If $N=2$ for a given $b>1$ we have from $(\ref{I-1-5*})$ that if
$$m<(p+1)\left(\frac{b-1}{bp}\right)^{\frac{p}{p+1}},
$$
then, for some $\gd>0$,
\bel{I-1-5-2-1}
-\Gd v\geq \gd v^{1+b(p-1)}.
\ee
The result follows from \rprop{zero} by choosing $b$ large enough.\qeda
\medskip

\subsubsection{ Proof of Theorem B'} 1- We assume that such a supersolution $u$ exists and we denote $u=e^v$, then
\bel{I-1-5-3-1}
-\Gd v\geq F(\abs{\nabla v}^2),
\ee
where
$$F(X)=X+e^{(p-1)v}+Me^{\frac{p-1}{p+1}v}X^{\frac{p}{p+1}}.$$
Clearly, if $M\geq 0$, then $F(X)\geq 0$ for any $X\geq 0$. Next we assume $M<0$, then
$$F(X)\geq F(X_0)=e^{(p-1)v}\left(1-p^p\left(\frac{\abs M}{p+1}\right)^{p+1}\right)=e^{(p-1)v}\left(1-\left(\frac{\abs M}{\gm^*(2)}\right)^{p+1}\right).
$$
Hence, if $\abs M\leq \gm^*(2)$,  $v$ is a positive superharmonic function in $\Gw$ which tends to infinity on the boundary. Such a function is larger than the harmonic function with boundary value $k>0$ for any $k$ (and taking the value $\displaystyle\min_{\abs x=R}v(x)$ for $R$ large enough if $\Gw$ is an exterior domain). Letting $k\to\infty$ we derive a contradiction.\smallskip

\nind 2- Let $R>0$ such that $\Gw^c\subset B_R$ and let $w$ be the solution of
\bel{I-1-5-1}\BA {lll}
-\Gd w-ae^{(p-1)w}=0\qquad&\text{in }B_R\cap \Gw\\\phantom{}
\!\!\!\!\!\!\displaystyle\lim_{\dist (x,\prt B_R)\to 0}w(x)=-\infty\\[2mm]\phantom{}
\!\!\displaystyle\lim_{\dist (x,\prt\Gw)\to 0}w(x)=\infty,
\EA\ee
with $a=1-\left(\frac{\abs M}{\gm^*(2)}\right)^{p+1}<0$, obtained by approximations. By the argument used in 1,
$$ae^{(p-1)w}\leq \abs{\nabla w}^2+e^{(p-1)w}-\abs Me^{\frac{p-1}{p+1}w}\abs{\nabla w}^{\frac{2p}{p+1}},
$$
hence
$$-\Gd w\leq \abs{\nabla w}^2+e^{(p-1)w}-\abs Me^{\frac{p-1}{p+1}w}\abs{\nabla w}^{\frac{2p}{p+1}}.
$$
Therefore $v=e^{w}$ is nonnegative and satisfies
\bel{I-1-5-2}\BA {lll}
-\Gd v-v^p+\abs M\abs{\nabla v}^{\frac{2p}{p+1}}\leq 0\qquad&\text{in }B_R\cap \Gw\\\phantom{-\Gd -v^p+\abs M\abs{\nabla v}^{\frac{2p}{p+1}}}
v=0\qquad&\text{on }\prt B_R\\\phantom{,\Gd--,, }
\!\!\!\!\displaystyle\lim_{\dist (x,\prt\Gw)\to 0}v(x)=\infty.
\EA\ee
Next we extend $v$ by zero in $B_R^c$ and denote by $\tilde v$ the new function. It is a nonnegative subsolution of $(\ref{I-0})$
which tends to $\infty$ on $\prt\Gw$. For constructing a supersolution we recall that if $M\leq -\gm^*(1)$ there exist two types of explicit solutions of
\bel{I-1-5-3}\BA {lll}
-u''=u^p+M\abs{u'}^{\frac{2p}{p+1}}
\EA\ee
defined on $\BBR$ by $U_{j,M}(t)=\infty$ for $t\leq 0$ and $U_{j,M}(t)=X_{j,M}t^{-\frac{2}{p-1}}$, j=1,2, for $t>0$ where $X_{1,M}$ and $X_{2,M}$ are respectively the smaller and the larger positive root of
\bel{I-1-5-4}\BA {lll}
X^{p-1}-\abs M\left({\myfrac{2}{p-1}}\right)^{\frac{2}{p+1}}X^{\frac{p-1}{p+1}}+\myfrac{2(p+1)}{(p-1)^2}=0.
\EA\ee

Since $\Gw^c$ is convex it is the intersection of all the closed half-spaces which contain it and we denote by $\CH_\Gw$ the family of such hyperplanes which are touching $\prt\Gw$. If $H\in \CH_\Gw$ let ${\bf n}_H$ be the normal direction to $H$, inward  with respect to $\Gw$, $\CH_+=\{x\in\BBR^N:\langle {\bf n}_H,x-{\bf n}_H\rangle>0\}$ and we define $U_H$ in the direction ${\bf n}_H$ by putting
$$U_H(x)=U_{2,M}(\langle {\bf n}_H,x-{\bf n}_H\rangle)=X_{2,M}\left(\langle {\bf n}_H,x-{\bf n}_H\rangle\right)^{-\frac{2}{p-1}}\quad\text{for all }\,x\in \CH_+.
$$
Hence
and set, for $x\in \Gw:=\cap_{H\in \CH_\Gw} \CH_+$,
\bel{I-1-5-5}\BA {lll}\displaystyle
u_\Gw(x)=\inf_{H\in \CH_\Gw}U_H(x).
\EA\ee
Then $u_\Gw$ is a nonnegative supersolution of $(\ref{I-0})$ in $\Gw$ and
$$u_\Gw(x)\geq X_{2,M}(\dist x,\Gw))^{-\frac{2}{p-1}}\qquad\forall x\in \Gw.
$$
Next $v_\Gw=\ln u_\Gw$ blows up on $\prt\Gw$, is finite on $\prt B_R$ and satisfies
\bel{I-1-5-6}-\Gd v_\Gw-ae^{(p-1)v_\Gw}\geq 0\qquad\text{in }B_R\cap \Gw.
\ee
By comparison with $w$ since $a<0$, $v_\Gw\geq w$. Hence $u_\Gw\geq v$ in $B_R\setminus\Gw^c$. Extending $v$ by zero as $\tilde v$ we obtain $u_\Gw\geq \tilde v$ in $\Gw^c$. Hence $u_\Gw$ is a supersolution in $\Gw^c$ where it dominates the subsolution
$\tilde v$. It follows by \cite[Theorem 1-4-6]{Vebook} that there exists a solution $u$ of $(\ref{I-0})$ satisfying $\tilde v\leq u\leq u_\Gw$, which ends the proof. \qeda

\mysection{The refined Bernstein method}
The method is a combination of the one used in the previous proofs. It is based upon the replacement of the unknown by setting first $u=v^{-\gb}$ as in \cite{GiSp} and \cite{BiVe} and the study of the equation satisfied by $\abs{\nabla v}$. However we do not use integral techniques.
Since $u$ is a positive solution of  $(\ref{I-0})$ in $B_R$, the function $v$ is well defined and satisfies
\bel{I-1-8}
-\Gd v+(1+\gb)\frac{\abs{\nabla v}^2}{v}+\frac{1}{\gb}v^{1-\gb(p-1)}+M\abs\gb^{q-2}\gb v^{(\gb+1)(1-q)}\abs{\nabla v}^q=0
\ee
in $B_R$. We set
$$z=\abs{\nabla v}^2\,,\;s=1-q-\gb(q-1)=(1-q)(\gb+1)\,,\;\gs=1-\gb(p-1),
$$
and derive
\bel{I-1-9}\Gd v=(1+\gb)\frac{z}{v}+\frac{1}{\gb}v^{\gs}+M\abs\gb^{q-2}\gb v^{s}z^{\frac{q}{2}}.
\ee
Combining Bochner's formula and Schwarz identity we have classically
$$\frac 12\Gd z\geq \frac 1N(\Gd v)^2+\langle\nabla\Gd v,\nabla v\rangle.
$$
We explicit the different terms
$$\BA {lll}
(\Gd v)^2=(1+\gb)^2\myfrac{z^2}{v^2}+M^2\gb^{2(q-1)}v^{2s}z^q+\myfrac{v^{2\gs}}{\gb^2}
+2M(1+\gb)\abs\gb^{q-2}\gb v^{s-1}z^{1+\frac q2}\\[2mm]
\phantom{---------------}
+\myfrac{2(1+\gb)}{\gb}v^{\gs-1}z+2M\abs\gb^{q-2}v^{s+\gs}z^{\frac q2},
\EA$$
$$\BA {lll}\nabla\Gd v=(1+\gb)\myfrac{\nabla z}{v}-\myfrac{(1+\gb)z}{v^2}\nabla v+\myfrac{\gs}{\gb}v^{\gs-1}\nabla v+Ms\abs\gb^{q-2}\gb v^{s-1}z^{\frac{q}{2}}\nabla v
\\[2mm]
\phantom{------------------}
+\myfrac{Mq}{2}\abs\gb^{q-2}\gb v^{s}z^{\frac{q}{2}-1}\nabla z,
\EA$$
$$\BA {lll}
\langle\nabla\Gd v,\nabla v\rangle=\left(\myfrac{1+\gb}{v}+\myfrac{Mq}{2}\abs\gb^{q-2}\gb v^{s}z^{\frac{q}{2}-1}\right)\langle\nabla z,\nabla v\rangle -\myfrac{(1+\gb)z^2}{v^2}+\myfrac{\gs}{\gb}v^{\gs-1}z\\[2mm]
\phantom{----------------------}+Ms\abs\gb^{q-2}\gb v^{s-1}z^{\frac{q}{2}+1}.
\EA
$$
Hence
\bel{I-1-10}\BA{lll}
-\myfrac{1}{2}\Gd z+\myfrac{1}{N}(\Gd v)^2+\left(\myfrac{1+\gb}{v}+\myfrac{Mq}{2}\abs\gb^{q-2}\gb v^{s}z^{\frac{q}{2}-1}\right)\langle\nabla z,\nabla v\rangle\\[2mm]
\phantom{-----------}
-\myfrac{(1+\gb)z^2}{v^2}+\myfrac{\gs}{\gb}v^{\gs-1}z+Ms\abs\gb^{q-2}\gb v^{s-1}z^{\frac{q}{2}+1}\leq 0.
\EA\ee

\subsection{ Proof of Theorem D}
We develop the term $(\Gd v)^2$ in $(\ref{I-1-10})$ and get
\bel{I-1-11}\BA {lll}
-\myfrac{1}{2}\Gd z +\left(\myfrac {(1+\gb)^2}{N}-(1+\gb)\right)\myfrac{z^2}{v^2}+\myfrac{M^2\gb^{2(q-1)}}{N}v^{2s}z^q
+M\left(s+\myfrac{2(1+\gb)}{N}\right)\abs\gb^{q-2}\gb v^{s-1}z^{1+\frac q2}\\[2mm]
+\myfrac{v^{2\gs}}{N\gb^2}+\left(\myfrac{1+\gb}{v}+\myfrac{Mq}{2}\abs\gb^{q-2}\gb v^{s}z^{\frac{q}{2}-1}\right)\langle\nabla z,\nabla v\rangle
+\myfrac{N\gs+2(1+\gb)}{N\gb}v^{\gs-1}z+\myfrac{2M\abs\gb^{q-2}}{N}v^{s+\gs}z^{\frac q2}\\[2mm]
\phantom{---------------------------------------}
\leq 0.
\EA\ee
Next we set $z=v^{-k}Y$ where $k$ is a real parameter. Then $\nabla z=-kv^{-k-1}Y\nabla v+v^{-k}\nabla Y$,
$$\langle\nabla z,\nabla v\rangle=-kv^{-k-1}Yz+v^{-k}\langle\nabla Y,\nabla v\rangle=-kv^{-2k-1}Y^2+v^{-k}\langle\nabla Y,\nabla v\rangle,
$$
$$\myfrac{\langle\nabla z,\nabla v\rangle}{v}=-kv^{-2k-2}Y^2+v^{-k-1}\langle\nabla Y,\nabla v\rangle,
$$
$$\BA {lll}Mv^sz^{\frac q2-1}\langle\nabla z,\nabla v\rangle
=-kMv^{s-\frac{qk}{2}-k-1}Y^{\frac q2+1}+Mv^{s-\frac{qk}{2}}Y^{\frac q2-1}\langle\nabla Y,\nabla v\rangle,
\EA$$
$$\BA {lll}
-\Gd z={\rm div}\left(kv^{-k-1}Y\nabla v-v^{-k}\nabla Y\right)
\\[1mm]\phantom{-\Gd z}
=kv^{-k-1}Y\Gd v-k(k+1)v^{-k-2}Yz+2kv^{-k-1}\langle\nabla Y,\nabla v\rangle-v^{-k}\Gd Y
\\[1mm]\phantom{-\Gd z}
=kv^{-k-1}Y\Gd v-k(k+1)v^{-2k-2}Y^2+2kv^{-k-1}\langle\nabla Y,\nabla v\rangle-v^{-k}\Gd Y.
\EA$$
From $(\ref{I-1-9})$
$$\Gd v=(1+\gb)v^{-k-1}Y+\frac{1}{\gb}v^{\gs}+M\abs\gb^{q-2}\gb v^{s-k\frac q2}Y^{\frac{q}{2}},$$
therefore
$$\BA {lll}
-\Gd z=k(\gb-k)v^{-2k-2}Y^2+\myfrac{k}{\gb}v^{\gs-k-1}Y+kM\abs\gb^{q-2}\gb v^{s-k\frac q2-k-1}Y^{\frac{q}{2}+1}\\
\phantom{-\Gd z--------------}
+2kv^{-k-1}\langle\nabla Y,\nabla v\rangle-v^{-k}\Gd Y.
\EA$$
Replacing $\langle\nabla z,\nabla v\rangle$ and $\Gd z$ given by the above expressions in $(\ref{I-1-11})$ and $z$ by $v^{-k}Y$, leads to
$$\BA {lll}
-\Gd Y+\left(\myfrac{k(\gb-k)}{2}+\myfrac{(1+\gb)^2}{N}-(k+1)(\gb+1)\right)v^{-k-2}Y^2+\myfrac{v^{2\gs+k}}{N\gb^2}
+\myfrac{M^2\gb^{2(q-1)}}{N}v^{2s+k-kq}Y^q\\[4mm]
\phantom{----}
+\left(\myfrac{k+\gb+1}{v}+\myfrac{Mq\abs\gb^{q-2}\gb}{2} v^{s+k-k\frac q2}Y^{\frac{q}{2}-1}\right)\langle\nabla Y,\nabla v\rangle
+ \myfrac{2M\abs\gb^{q-2}}{N}v^{s+\gs+k-k\frac q2}Y^{\frac q2}\\[4mm]
+\left(s+\myfrac{2(1+\gb)}{N}-\myfrac{k(q-1)}{2}\right)M\abs\gb^{q-2}\gb v^{s-k\frac q2-1}Y^{1+\frac q2}
+\myfrac{1}{\gb}\left(\myfrac{k}{2}+\gs+\myfrac{2(1+\gb)}{N}\right)v^{\gs-1}Y\leq 0.
\EA$$
For $\ge_1,\ge_2>0$,
$$\myfrac{1}{v}\abs{\langle\nabla Y,\nabla v\rangle}\leq \ge_1 v^{-k-2}Y^2+\myfrac{1}{4\ge_1}\myfrac{\abs{\nabla Y}^2}{Y},
$$
$$v^{s+k-k\frac q2}Y^{\frac{q}{2}-1}\abs{\langle\nabla Y,\nabla v\rangle}\leq \ge_2v^{2s-kq+k}Y^{q}+\myfrac{1}{4\ge_2}\myfrac{\abs{\nabla Y}^2}{Y}.
$$
Hence
\bel{I-1-12}\BA {lll}
-\Gd Y+\myfrac{v^{2\gs+k}}{N\gb^2}+ \myfrac{2M\abs\gb^{q-2}}{N}v^{s+\gs+k-k\frac q2}Y^{\frac q2}
+\left(\myfrac{M^2\gb^{2(q-1)}}{N}-\myfrac{Mq\ge_2\abs\gb^{q-1}}{2}\right)v^{2s+k-kq}Y^q\\[4mm]\phantom{-\Gd Y}
+\left(\myfrac{k(\gb-k)}{2}+\myfrac{(1+\gb)^2}{N}-(k+1)(\gb+1)-\abs{k+\gb+1}\ge_1\right)v^{-k-2}Y^2\\[4mm]\phantom{-\Gd Y}
+\myfrac{1}{\gb}\left(\myfrac{k}{2}+\gs+\myfrac{2(1+\gb)}{N}\right)v^{\gs-1}Y+\left(s+\myfrac{2(1+\gb)}{N}-\myfrac{k(q-1)}{2}\right)M\abs\gb^{q-2}\gb v^{s-k\frac q2-1}Y^{1+\frac q2}\\[4mm]\phantom{-\Gd Y}
\leq \left(\myfrac{\abs{k+\gb+1}}{\ge_1}+\myfrac{Mq\abs\gb^{q-1}}{2\ge_2}\right)\myfrac{\abs{\nabla Y}^2}{4Y}.
\EA\ee
We first choose $\ge_2=\frac{M\abs\gb^{q-1}}{qN}$, then
\bel{I-1-12x}\BA {lll}
-\Gd Y+\myfrac{v^{2\gs+k}}{N\gb^2}+\left(\myfrac{k(\gb-k)}{2}+\myfrac{(1+\gb)^2}{N}-(k+1)(\gb+1)-\abs{k+\gb+1}\ge_1\right)v^{-k-2}Y^2
\\\phantom{-\Gd Y}
+\myfrac{1}{\gb}\left(\myfrac{k}{2}+\gs+\myfrac{2(1+\gb)}{N}\right)v^{\gs-1}Y+
\myfrac{M^2\gb^{2(q-1)}}{2N}v^{2s+k-kq}Y^q+
\myfrac{2M\abs\gb^{q-2}}{N}v^{s+\gs+k-k\frac q2}Y^{\frac q2}
\\[4mm]\phantom{-\Gd Y}
+\left(s+\myfrac{2(1+\gb)}{N}-\myfrac{k(q-1)}{2}\right)M\abs\gb^{q-2}\gb v^{s-k\frac q2-1}Y^{1+\frac q2}\\[4mm]\phantom{-\Gd Y}
\leq \left(\myfrac{\abs{k+\gb+1}}{\ge_1}+\myfrac{Nq^2}{2}\right)\myfrac{\abs{\nabla Y}^2}{4Y}.
\EA\ee
In order to show the sign of the terms on the left in $(\ref{I-1-12})$, we separate the terms containing the coefficient $M$ from the ones which do not contain it. Indeed these last terms are associated to the mere Lane-Emden equation $(\ref{I-1})$ which is treated, as a particular case, in \cite[Theorem B]{BiGaVe1} where the exponents therein are $q=0$, and $p\in \left(1,\frac{N+3}{N-1}\right)$. We set
\bel{I-1-12^*}\BA {lll}
H_{\ge_1,1}=\myfrac{v^{2\gs+k}}{N\gb^2}+\left(\myfrac{k(\gb-k)}{2}+\myfrac{(1+\gb)^2}{N}-(k+1)(\gb+1)-\abs{k+\gb+1}\ge_1\right)v^{-k-2}Y^2\\[2mm]
\phantom{--------------}+\myfrac{1}{\gb}\left(\myfrac{k}{2}+\gs+\myfrac{2(1+\gb)}{N}\right)v^{\gs-1}Y\\[2mm]
\phantom{H_{\ge_0,1}}
=v^{2\gs+k}\tilde H_{\ge_1,1}(v^{-1-k-\gs}Y),
\EA\ee
where
\bel{I-1-13}\BA {lll}
\tilde H_{\ge_1,1}(t)=\left(\myfrac{k(\gb-k)}{2}+\myfrac{(1+\gb)^2}{N}-(k+1)(\gb+1)-\abs{k+\gb+1}\ge_1\right)t^2\\[4mm]
\phantom{\tilde H_{\ge_1,1}(t)-----------}
+\myfrac{1}{\gb}\left(\myfrac{k}{2}+\gs+\myfrac{2(1+\gb)}{N}\right)t+\myfrac{1}{N\gb^2},
\EA\ee
and
\bel{I-1-14}\BA {lll}
H_{M,2}=\myfrac{M^2\gb^{2(q-1)}}{2N}v^{2s+k-kq}Y^q+\myfrac{2M\abs\gb^{q-2}}{N}v^{s+\gs+k-k\frac q2}Y^{\frac q2}\\[4mm]
\phantom{------}
+\left(s+\myfrac{2(1+\gb)}{N}-\myfrac{k(q-1)}{2}\right)M\abs\gb^{q-2}\gb v^{s-k\frac q2-1}Y^{1+\frac q2}.
\EA\ee
Then
$$\BA {lll}
-\Gd Y+v^{2\gs+k}\tilde H_{\ge_1,1}(v^{-1-k-\gs}Y)+H_{M,2}\leq \left(\myfrac{\abs{k+\gb+1}}{\ge_1}+\myfrac{Nq^2}{2}\right)\myfrac{\abs{\nabla Y}^2}{4Y}.
\EA$$
The sign of $\tilde H_{\ge_1,1}$ depends on its discriminant $\CD_{\ge_1}$ which is a polynomial in its coefficients. Then if for $\ge_1=0$ this discriminant is negative $\CD_{0}$ is negative, the discriminant $\CD_{\ge_1}$ of  $\tilde H_{\ge_1,1}$ shares this property for $\ge_1>0$ small enough and therefore $H_{\ge_1,1}$ is positive.
The proof is similar as the one of  \cite[Theorem B]{BiGaVe1} in case
(i) but for the sake of completeness we recall the main steps. Firstly
$$\BA{lll}
\CD'_0:=N^2\gb^2\CD_0=\left(\myfrac{Nk}{2}+\gs N+2(1+\gb)\right)^2-4\left(\myfrac{Nk(\gb-k)}{2}+(1+\gb)^2-N(k+1)(\gb+1)\right).
\EA$$
Then
$$\CD'_0=\left(\myfrac {N(p-1)}4-1\right)(2\gs+k)^2+2(p-1)(2\gs+k)+\tilde L
$$
where $\tilde L=(p-1)k^2+p(\gl+2)^2>0$. Put
$$S=\myfrac{2\gs+k}{k+2}=1-\myfrac{2\gb(p-1)}{k+2}\,\text{ and }\; \CT(S)=\left(\myfrac{(N-1)(p-1)}{4}-1\right)
S^2+(p-1)S+p.
$$
After some computations we get, if $k\neq -2$,
\bel{I-1-14a}\BA {llll}
\CD'_1:=\myfrac{(p-1)\CD'_0}{(k+2)^2}=(p-1)\left(\myfrac{k}{k+2}-\myfrac{S}{2}\right)^2 +\CT(S).
\EA\ee
We choose $S>2$ such that $\frac{k}{k+2}-\frac{S}{2}=0$, hence $\gb=\frac{2-k}{2(p-1)}$. If $p<\frac{N+3}{N-1}$ the coefficient of $S^2$ in $\CT(S)$ is negative. Hence $\CT(S)<0$ provided $S$ is large enough which is satisfied if $k<-2$ with $\abs{k+2}$ small enough. We infer from this that $\gb>0$, $\CD_0<0$ and $\tilde H_{\ge_1,1}>0$ if $\ge_1$ is small enough. In particular $\tilde H_{\ge_1,1}(t)\geq c_6 (t^2+1)$
for some $c_6=c_6(N,p,q)>0$, which means
\bel{I-1-15}
v^{2\gs+k}\tilde H_{\ge_1,1}(v^{-1-k-\gs}Y)\geq c_6\left(v^{-k-2}Y^2+v^{2\gs+k}\right).
\ee

\nind Secondly the positivity of $H_{M,2}$ is ensured, as $\gb$ and $M$ are positive, by the positivity of
$$\CA:=s+\myfrac{2(1+\gb)}{N}-\myfrac{k(q-1)}{2}.
$$
Replacing $s$ by its value, we obtain, since $1<q<\frac{N+2}N$ and $\gb+\frac{2+k}{2}>0$, which can be assume by taking $\abs{k+2}$ small enough,
$$\CA=2\myfrac{1+\gb}{N}-(q-1)\left(\gb+1+\myfrac{k}{2}\right)>-\myfrac{k}{N}
$$
Then we deduce that
\bel{I-1-16}\BA {lll}
-\Gd Y +c_6\left(v^{-k-2}Y^2+v^{2\gs+k}\right)\leq c_7\myfrac{\abs{\nabla Y}^2}{Y},
\EA\ee
and $c_7=c_7(N,p,q)>0$ is independent of $M$. Since  $S=1-\frac{2\gb (p-1)}{k+2}=1-\frac{2-k}{k+2}=\frac{2k}{k+2}>0$, we have
\bel{I-1-17}\BA {lll}
2Y^{\frac{2S}{S+1}}=2\left(\myfrac{Y^2}{v^{k+2}}\right)^{\frac{S}{S+1}}v^{\frac{(k+2)S}{S+1}}\leq \myfrac{Y^2}{v^{k+2}}+v^{(k+2)S}=\myfrac{Y^2}{v^{k+2}}+v^{2\gs+k}.
\EA\ee
From this we infer the inequality
\bel{I-1-18}\BA {lll}
-\Gd Y +2c_6Y^{\frac{2S}{S+1}}\leq c_7\myfrac{\abs{\nabla Y}^2}{Y}.
\EA\ee
Then we derive from \rlemma{L2.2}  that in the ball $B_R$ there holds
\bel{I-1-19}\BA {lll}
Y(0)\leq c_8R^{-\frac{2(S+1)}{S-1}}=c_8R^{-2+\frac{2(k+2)}{\gb(p-1)}}.
\EA\ee
From this it follows
\bel{I-1-20}\BA {lll}
\abs{\nabla u^{-\frac{2+k}{2\gb}}(0)}\leq \myfrac{\abs{k+2}}{2}\sqrt{c_8\,}R^{-1+\frac{k+2}{\gb(p-1)}}.
\EA\ee
Setting $a=-\frac{k+2}{2\gb}>0$ we get that for any domain $\Gw\subset\BBR^N$ any positive solution in $\Gw$ satisfies
\bel{I-1-21}\BA {lll}
\abs{\nabla u^{a}(x)}\leq \myfrac{|k+2|}{2}\sqrt{c_8\,}\left(\dist(x,\prt\Gw)\right)^{-1-\frac{2a}{p-1}}\qquad\text{for all }\, x\in\Gw.
\EA\ee
The non existence of any positive of $(\ref{I-0})$ solution in $\BBR^N$ follows classically. \phantom{--}\phantom{--}\qeda

\bcor{B2} Let $\Gw$ be a smooth domain in $\BBR^N$, $N \geq  2$ with a bounded boundary, $1<p<\frac{N+3}{N-1}$, $1<q<\frac{N+2}{N}$ and $M>0$. If $u$ is a positive solution of $(\ref{I-0})$ in $\Gw$ there exists $d_0$ depending on $\Gw$ and $c_9 = c_9(N,p,q) > 0$ such that
\bel{I-1-21x}\BA {lll}\displaystyle
u(x)\leq c_9\left(\left(\dist(x,\prt\Gw)\right)^{-\frac{2}{p-1}}+\max_{\dist(z,\prt\Gw)=d_0}u(z)\right)\qquad\text{for all }\, x\in\Gw.
\EA\ee
\es

\nind\Proof It is similar to the one of \cite[Corollary B-2]{BiGaVe1}.\qeda
\mysection{The integral method}
\subsection{Preliminary inequalities}
We recall the next  inequality \cite[Lemma 3.1]{BiRa}.
\blemma{BR1} Let $\Gw\subset\BBR^N$ be a domain. Then for any positive $u\in C^2(\Gw)$, any nonnegative $\eta\in C^\infty_0(\Gw)$ and any real numbers $m$ and $d$ such that $d\neq m+2$, the following inequality holds
\bel{I-2-1}\BA {lll}
A\myint{\Gw}{}\eta u^{m-2}\abs{\nabla u}^4 dx-\myfrac{N-1}{N}\myint{\Gw}{}\eta u^{m}(\Gd u)^2 dx
-B\myint{\Gw}{}\eta u^{m-1}\abs{\nabla u}^2 \Gd u dx\leq R,
\EA\ee
where
$$\BA {lll}
A=\myfrac{1}{4N}\left(2(N-m)d-(N-1)(m^2+d^2)\right)\,,\; B=\myfrac{1}{2N}\left(2(N-1)m+(N+2)d\right),
\EA$$
and
$$R=\myfrac{m+d}{2}\myint{\Gw}{} u^{m-1}\abs{\nabla u}^2\langle\nabla u,\nabla\eta\rangle dx
+\myint{\Gw}{}u^{m}\Gd u \langle\nabla u,\nabla\eta\rangle dx+
\myfrac{1}{2}\myint{\Gw}{}u^{m}\abs{\nabla u}^2\Gd\eta dx.
$$
\es
It is noticeable that $d$ is a free parameter which plays a role only in the coefficients of the integral terms. The following technical result is useful to deal with the multi-parameter constraints problems which occur in our construction. It was first used in \cite{BiVe} under a simpler form and extended in \cite[Lemma 3.4]{BiRa}.
\blemma{BR2} For any $N\in\BBN$, $N\geq 3$ and $1<p<\frac{N+2}{N-2}$ there exist real numbers $m$ and $d$ verifying
\bel{I-2-2}\BA {lll}
(i)\qquad\qquad\qquad\qquad &d\neq m+2,\qquad\qquad\qquad\qquad\qquad\qquad\qquad\qquad\qquad\phantom{-}\\[4mm]
(ii) \qquad\qquad & \myfrac{2(N-1)p}{N+2}<d,\\[4mm]
(iii) \qquad\qquad &\max\left\{-2,1-p,\myfrac{(N-4)p-N}{2}\right\}<m\leq 0,\\[4mm]
(iv)\qquad\qquad & 2(N-m)d-(N-1)(m^2+d^2)>0.
\EA\ee
\es
\subsection{Proof of Theorem E}

{\it Step 1: The integral estimates}. Let $\eta\in C^{\infty}_0(\Gw)$, $\eta\geq 0$. We apply \rlemma{BR1} to a positive solution $u\in C^2(\Gw)$ of $(\ref{I-0})$, firstly with $q>1$ and then with $q=\frac{2p}{p+1}$.
\bel{I-2-3}\BA {lll}
A\myint{\Gw}{}\eta u^{m-2}\abs{\nabla u}^4 dx-\myfrac{N-1}{N}\myint{\Gw}{}\eta\left(u^{m+2p}+2Mu^{m+p}\abs{\nabla u}^q+M^2u^{m}\abs{\nabla u}^{2q}\right) dx\\[4mm]\phantom{---------------------}
-B\myint{\Gw}{}\eta u^{m-1}\abs{\nabla u}^2 \Gd u dx\leq R.
\EA\ee
We multiply $(\ref{I-0})$ by $\eta u^{m+p}$ and integrate over $\Gw$. Then
$$\BA {lll}
\myint{\Gw}{}\eta \left(u^{m+2p}+Mu^{m+p}\abs{\nabla u}^q \right)dx=-\myint{\Gw}{}\eta u^{m+p}\Gd udx\\[4mm]
\phantom{\myint{\Gw}{}\eta \left(u^{m+2p}+Mu^{m+p}\abs{\nabla u}^q \right)dx}
=\myint{\Gw}{}u^{m+p}\langle\nabla u,\nabla\eta\rangle dx+(m+p)\myint{\Gw}{}\eta u^{m+p-1}\abs{\nabla u}^2 dx.
\EA$$
We set
$$\BA {lll}
F=\myint{\Gw}{}\eta u^{m-2}\abs{\nabla u}^4 dx\,,\;P=\myint{\Gw}{}\eta u^{m-1}\abs{\nabla u}^{q+2} dx\,,\;V=\myint{\Gw}{}\eta u^{m+2p} dx,\\[4mm]
T=\myint{\Gw}{}\eta u^{m+p-1}\abs{\nabla u}^{2} dx\,,\;W=\myint{\Gw}{}\eta u^{m+p}\abs{\nabla u}^{q} dx\,,\;U=\myint{\Gw}{}\eta u^{m} \abs{\nabla u}^{2q}dx,\\[4mm]
S=\myint{\Gw}{} u^{m+p}\langle \nabla u,\nabla\eta\rangle dx,
\EA$$
so that there holds
\bel{I-2-4}\BA {lll}
AF-\myfrac{N-1}{N}\left(V+2MW+M^2U\right)+BT+BMP\leq R,
\EA\ee
and
\bel{I-2-5}\BA {lll}
V+MW=(m+p)T+S.
\EA\ee
Eliminating $V$ between $(\ref{I-2-4})$ and $(\ref{I-2-5})$, we get
\bel{I-2-6}\BA {lll}
AF+B_0T+M\left(BP-\myfrac{N-1}{N}W-\myfrac{N-1}{N}MU\right)\leq R-\myfrac{N-1}{N}S,
\EA\ee
where
$$B_0=B-\myfrac{N-1}{N}(m+p)=\myfrac{N+2}{2N}d-\myfrac{N-1}{N}p.
$$
Also
$$2P=2\myint{\Gw}{}\eta u^{m}\myfrac{\abs{\nabla u}^{2}}{u}\abs{\nabla u}^{q} dx\leq
\myint{\Gw}{}\eta u^m\left(\myfrac{\abs{\nabla u}^{4}}{u^2}+\abs{\nabla u}^{2q}\right) dx=F+U.
$$
We fix now $q=\frac{2p}{p+1}$, then
\bel{I-2-8}\BA {lll}
U=\myint{\Gw}{}\eta u^m\abs{\nabla u}^{2q}dx=\myint{\Gw}{}\eta u^m\left(\myfrac{\abs{\nabla u}}{\sqrt u}\right)^{4(q-1)}u^{2(q-1)}\abs{\nabla u}^{4-2q}dx\\[4mm]\phantom{U}
\leq \myfrac{p-1}{p+1}\myint{\Gw}{}\eta u^{m-2}\abs{\nabla u}^{4}dx+\myfrac{2}{p+1}\myint{\Gw}{}\eta u^{m+p-1}\abs{\nabla u}^{2}dx\\[4mm]\phantom{U}
\leq \myfrac{p-1}{p+1}F+\myfrac{2}{p+1}T,
\EA\ee
hence
\bel{I-2-7}\BA {lll}
P\leq \myfrac 12F+\myfrac 12U\leq  \myfrac {p}{p+1}F+\myfrac {1}{p+1}T
\EA\ee
and
\bel{I-2-9}\BA {lll}
2W=2\myint{\Gw}{}\eta u^{m+p}\abs{\nabla u}^{q}dx\leq\myint{\Gw}{}\eta u^{m+2p}dx+\myint{\Gw}{}\eta u^{m}\abs{\nabla u}^{2q}dx=V+U\\[4mm]
\phantom{W}
\leq U+(m+p)T+S-MW.
\EA\ee
Next we assume that $|M|\leq 1$. From $(\ref{I-2-8})$, $(\ref{I-2-9})$, it follows that
\bel{I-2-10}\BA {lll}
W\leq U+(m+p)T+S\leq F+(m+p+1)T+S.
\EA\ee
From now we fix $m$ and $d$ according \rlemma {BR2}. Therefore $A>0$ by $(\ref{I-2-2})$-(iv) and $B>0$ by combining $(\ref{I-2-2})$-(ii) and $(\ref{I-2-2})$-(iii). Furthermore $B_0>0$ by $(\ref{I-2-2})$-(ii).
Hence, from $(\ref{I-2-8})$, $(\ref{I-2-7})$ and $(\ref{I-2-10})$ we derive, since $\frac {N-1}{N}<1$ and $m\leq 0$ from   $(\ref{I-2-2})$-(ii)
$$\BA {lll}
\left | BP-\myfrac{N-1}{N}W-\myfrac{N-1}{N}MU\right |\leq  B\left(F+T\right)+F+(p+1)T+S+F+T,\\[4mm]
\phantom{  \left | BP-\myfrac{N-1}{N}W-\myfrac{N-1}{N}MU\right |}
\leq \left(B+2\right)F+\left(B+p+2\right)T+S.
\EA$$
Plugging these estimates into  $(\ref{I-2-6})$ we infer
\bel{I-2-11}\BA {lll}
AF+B_0T-|M|\left(\left(B+2\right)F+\left(B+p+2\right)T+S^{\phantom {e^d}}\!\!\!\!\right)\leq R-\myfrac{N-1}{N}S.
\EA\ee
Since $A$ and $B_0$ are positive, there exists $\gm_1\in (0,1)$ such that for any $\abs M<\gm_1$,
$$A_1:=A-|M|\left(B+2\right)>\frac A2\quad\text{and }\; B_1:=B_0-\abs M\left(B+p+2\right)>\frac{B_0}{2}.$$
Set $A_2=\min\{A_1,B_1\}$, then, and whatever is the sign of $S$,
$$A_2(F+T)\leq \abs R+\abs S.
$$
Using $(\ref{I-2-8})$ and $(\ref{I-2-7})$ we have
\bel{I-2-12}\BA {lll}
A_2(U+P)\leq  2A_2(F+T)\leq2(|R|+|S|).
\EA\ee
In the sequel we denote by $c_j$ some positive constants depending on $N$ and $p$. Then
\bel{I-2-13}\BA {lll}
U+P+F+T+W\leq c_{1}(|R|+|S|).
\EA\ee
On the other hand, we have
$$\abs R\leq
c_{2}\myint{\Gw}{} \left(u^{m-1}\abs{\nabla u}^3\abs{\nabla\eta} +u^{m+p}\abs{\nabla u}\abs{\nabla\eta} +u^{m}\abs{\nabla u}^{q+1}\abs{\nabla\eta} +u^{m}\abs{\nabla u}^2\abs{\Gd\eta} \right)dx.
$$
Since
$$\abs{\nabla u}^q=\left(\myfrac{\abs{\nabla u}}{\sqrt u}\right)^qu^{\frac q2}\leq \myfrac{\abs{\nabla u}^2}{ u}+
u^{\frac q{2-q}}=\myfrac{\abs{\nabla u}^2}{ u}+
u^{p},
$$
we deduce
$$\BA {lll}
\myint{\Gw}{} u^m|\nabla u|^{q+1}|\nabla \eta|dx\leq \myint{\Gw}{} u^{m-1}|\nabla u|^{3}|\nabla \eta|dx
+\myint{\Gw}{} u^{m+p}|\nabla u||\nabla \eta|dx.
\EA$$
Thus we derive from $(\ref{I-2-13})$
\bel{I-2-14}\BA {lll}
U+P+F+T+W\leq 2c_{3}\left(\myint{\Gw}{} u^{m-1}|\nabla u|^{3}|\nabla \eta|dx+\myint{\Gw}{} u^{m+p}|\nabla u||\nabla \eta|dx\right.\\[4mm]
\phantom{U+P+F+T+W---+\myint{\Gw}{} u^{m}\abs{\nabla u}^{q+1}\abs{\nabla\eta}dx}
\left.  +\myint{\Gw}{} u^{m}\abs{\nabla u}^2\abs{\Gd\eta} dx
\right).
\EA\ee
From this point we can use the method developed in \cite[p 599]{BiVe} for proving the Harnack inequality satisfied by positive solutions of
$(\ref{I-1})$ in $\Gw$. We set $\eta=\xi^\gl$ with  $\xi\in C^\infty_0(\Gw)$ with value in $[0,1]$ and $\gl>4$. For  $\ge\in (0,1)$ we have by the H\"older-Young inequality
\bel{I-2-15}\BA {lll}\myint{\Gw}{} u^{m-1}|\nabla u|^{3}|\nabla \xi^\gl|dx\leq\myfrac{\ge}{4c_{3}}\myint{\Gw}{} u^{m-2}|\nabla u|^{4}\xi^\gl dx
+C(\ge,c_3)\myint{\Gw}{} u^{m+2}|\nabla \xi|^{4}\xi^{\gl-4}dx,
\EA\ee
\bel{I-2-16}\BA {lll}\myint{\Gw}{} u^{m+p}|\nabla u||\nabla \xi^p|dx\leq\myfrac{\ge}{4c_{3}}\myint{\Gw}{} u^{m+p-1}|\nabla u|^{2}\xi^pdx
+C(\ge,c_3)\myint{\Gw}{} u^{m+p+1}|\nabla \xi|^{2}\xi^{\gl-2}dx,
\EA\ee
and
\bel{I-2-17}\BA {lll}\myint{\Gw}{} u^{m}|\nabla u|^2|\Gd \xi^p|dx\leq\myfrac{\ge}{4c_{3}}\myint{\Gw}{} u^{m-2}|\nabla u|^{4}\xi^pdx
+C(\ge,c_3)\myint{\Gw}{} u^{m+2}\left(|\nabla \xi|^{4}+\abs{\Gd \xi}^2\right)\xi^{\gl-4} dx.
\EA\ee
Hence
\bel{I-2-18}\BA {lll}
U+P+F+T+W\leq c_{4}  \left(\myint{\Gw}{} u^{m+2}\left(|\nabla \xi|^{4}+\abs{\Gd \xi}^2\xi^{2}\right)\xi^{\gl-4} dx
+\myint{\Gw}{} u^{m+p+1}|\nabla \xi|^{2}\xi^{\gl-2}dx \right).
\EA\ee
Let us denote by $c_{4}X$ the right-hand side of $(\ref{I-2-18})$. Combining  $(\ref{I-2-5})$, $(\ref{I-2-16})$ and $(\ref{I-2-18})$ we also get
\bel{I-2-19}\BA {lll}
S:= \myint{\Gw}{} u^{m+p}|\nabla u||\nabla \xi^p|dx\leq c_{5}X\Longrightarrow V:=\myint{\Gw}{} u^{m+2p}\xi^pdx\leq c_{6}X,
\EA\ee
and we finally obtain
\bel{I-2-20}\BA {lll}
U+V+P+F+S+T+W\leq c_{7} X.
\EA\ee
Finally we estimate the different terms in $X$, using that $m+p>0$ from $(\ref{I-2-2})$-(iii).  For $\ge>0$
\bel{I-2-21}\BA {lll}
\myint{\Gw}{} u^{m+2}\left(|\nabla \xi|^{4}+\abs{\Gd \xi}^2\xi^{2}\right)\xi^{\gl-4} dx\leq \ge \myint{\Gw}{} u^{m+2p}\xi^{\gl} dx\\[4mm]
\phantom{----------------}
+C(\ge,c_7)\myint{\Gw}{}\xi^{\gl-2\frac{m+2p}{p-1}}\left(|\nabla \xi|^{4}+\abs{\Gd \xi}^2\right)^{\frac{m+2p}{2(p-1)}}dx,
\EA\ee
and
\bel{I-2-22}\BA {lll}
\myint{\Gw}{} u^{m+p+1}|\nabla \xi|^{2}\xi^{\gl-2}dx\leq \ge \myint{\Gw}{} u^{m+2p}\xi^{\gl} dx
+C(\ge,c_7)\myint{\Gw}{}\xi^{\gl-2\frac{m+2p}{p-1}} |\nabla \xi|^{\frac{2(m+2p)}{p-1}}dx.
\EA\ee
At end we obtain
\bel{I-2-23}\BA {lll}
U+V+P+F+S+T+W\leq c_{8} \myint{\Gw}{}\xi^{\gl-2\frac{m+2p}{p-1}}\left(|\nabla \xi|^{4}+\abs{\Gd \xi}^2\right)^{\frac{m+2p}{2(p-1)}}dx.
\EA\ee

\nind {\it Step 2: The Harnack inequality}. We suppose that $\Gw=B_R\setminus\{0\}:=B_R^*$, fix $y\in B_{\frac R2}^*$, set $r=|y|$, then
$B_r(y)\subset B_R^*$. Let $\xi\in C^\infty_0(B_r(y))$ such that $0\leq \xi\leq 1$, $\xi=1$ in $B_{\frac r2}(y)$, $\abs{\nabla\xi}\leq cr^{-1}$ and $\abs{\Gd\xi}\leq cr^{-2}$. We choose $\gl>\max\left\{4,\frac{m+2p}{p+1}\right\}$, then
$$\myint{B_r(y)}{}\xi^{\gl-2\frac{m+2p}{p-1}}\left(|\nabla \xi|^{4}+\abs{\Gd \xi}^2\right)^{\frac{m+2p}{2(p-1)}}dx
\leq c_{9}r^{N-\frac{2(m+2p)}{p-1}},
$$
and hence
\bel{I-2-24}\BA {lll}
\myint{B_{\frac r2}(y)}{}u^{m+2p}dx\leq V\leq c_{10}r^{N-\frac{2(m+2p)}{p-1}}.
\EA\ee
We write $(\ref{I-0})$ under the form
\bel{I-2-25}\BA {lll}
\Gd u+D(x)u+M\langle G(x).\nabla u\rangle=0,
\EA\ee
with
$$D(x)=u^{p-1}\quad\text{and }\;\,G(x)=|\nabla u|^{-\frac{2}{p+1}}\nabla u.
$$
Set $\gs=\frac{m+2p}{p-1}$, then $\gs>\frac{N}{2}$ by $(\ref{I-2-2})$-(iii) and
\bel{I-2-26}\BA {lll}
\myint{B_{\frac r2}(y)}{}D^{\gs}(x)dx\leq V\leq c_{10}r^{N-\frac{2(m+2p)}{p-1}}=c_{10}r^{N-2\gs}.
\EA\ee
Next we estimate $G$. For $\gt,\gw,\gg>0$ and $\gth>1$, we have with $\gth'=\frac\gth{\gth-1}$,
$$\BA {lll}
\abs{\nabla u}^{(q-1)\gt}=u^{\gw}\abs{\nabla u}^{\gg}u^{-\gw}\abs{\nabla u}^{(q-1)\gt-\gg}\leq u^{\gw\gth'}\abs{\nabla u}^{\gg\gth}+u^{-\gw\gth}\abs{\nabla u}^{((q-1)\gt-\gg)\gth'}.
\EA$$
We fix
$$\gt=2\myfrac{2p+m}{p-1}=2\gs\,,\; \gw=\myfrac{(2-m)(p+m-1)}{p+1}\,\text{ and }\; \gth=\myfrac{p+1}{2-m}.
$$
Then $\gw>0$ and $\gth>1$ from $(\ref{I-2-2})$-(iii), $\gw>0$. Then $u^{\gw\gth'}\abs{\nabla u}^{\gg\gth}=u^{p+m-1}\abs{\nabla u}^{2}$ and
$u^{-\gw\gth}\abs{\nabla u}^{((q-1)\gt-\gg)\gth'}=u^{m-2}\abs{\nabla u}^{4}$, thus
$$\myint{B_{\frac r2}(y)}{}\abs{\nabla u}^{(q-1)\gt}dx\leq F+T\leq c_{11} \myint{\Gw}{}\xi^{\gl-2\frac{m+2p}{p-1}}\left(|\nabla \xi|^{4}+\abs{\Gd \xi}^2\xi^{2}\right)^{\frac{m+2p}{2(p-1)}}dx.
$$
This implies
\bel{I-2-27}\BA {lll}
\myint{B_{\frac r2}(y)}{}G^{\gt}(x)dx\leq c_{12}r^{N-\gt},
\EA\ee
with $\gt>N$.   Using the results of \cite[Section 5]{Tru}, we infer that a Harnack inequality, uniform with respect to $r$, is satisfied. Hence there  exists $c_{13}>0$ depending on $N,p$ such that for any $r\in (0,\frac R2]$ and $y$ such that $\abs y=r$ there holds
\bel{I-2-28}\BA {lll}\displaystyle
\max_{z\in B_{\frac r2}(y)}u(z)\leq c_{13}\min_{z\in B_{\frac r2}(y)}u(z)\quad\forall 0<r\leq \tfrac R2\;\;\forall y\,\text{ s.t. }\,|y|=r,
\EA\ee
which implies
\bel{I-2-29}\BA {lll}\displaystyle
u(x)\leq c_{14} u(x')\quad\forall x,x'\in\BBR^N\;\;\text{ s.t. }\,|x|=|x'|\leq \frac R2,
\EA\ee
and actually $c_{14}=c_{13}^{7}$ by a simple geometric construction.  By  $(\ref{I-2-24})$
$$r^N\gw_Nr^N\left(\min_{z\in B_{\frac r2}(y)}u(z)\right)^{m+2p}\leq 4^Nc_{10}r^{N-\frac{2(m+2p)}{p-1}},
$$
where $\gw_N$ is the volume of the unit N-ball. This implies
\bel{I-2-30}\BA {lll}\displaystyle
u(x)\leq c_{14}\abs x^{-\frac{2}{p-1}}\qquad\forall x\in B_{\frac R2}^*.
\EA\ee
The proof follows.\qeda\medskip

\nind{\Remark} Using standard rescaling techniques (see e.g. \cite[Lemma 3.3.2]{Vebook})  the gradient estimate holds
\bel{I-2-31}\BA {lll}\displaystyle
\abs{\nabla u(x)}\leq c_{15}\abs x^{-\frac{p+1}{p-1}}\qquad\forall x\in B_{\frac R3}^*.
\EA\ee
And the next estimate for a solution $u$ in a domain $\Gw$ satisfying the interior sphere condition with radius $R$ is valid
\bel{I-2-32}\BA {lll}\displaystyle
u(x)\leq c_{14}\left(\dist( x,\prt\Gw)\right)^{-\frac{2}{p-1}}\quad\forall x\in \Gw\;\text{ s.t.}\; \dist( x,\prt\Gw)\leq \frac R2.
\EA\ee
\mysection {Radial ground states}
We recall that if $q\neq\frac{2p}{p+1}$ and $M\neq 0$, $(\ref{I-0})$ can be reduced to the case $M=\pm 1$ by using the transformation $(\ref{I-9})$.
Since any ground state $u$ of $(\ref{I-0})$ radial with respect to $0$ is decreasing (this is classical and straightforward), it achieves its maximum at $0$ and the following equivalence holds if $v$ is defined by $(\ref{I-9})$
\bel{X-4}\BA {lll}\displaystyle
-u''-\myfrac{N-1}{r}u'=|u|^{p-1}u+M\abs{u_r}^q&\quad\text{s.t.  }\, &\max u=u(0)=1\\
\Longleftrightarrow \\
-v''-\myfrac{N-1}{r}v'=|v|^{p-1}v\pm\abs{v_r}^q&\quad\text{s.t. }\, &\max v=v(0)=|M|^{\frac{2}{(p+1)q-2p}}.
\EA\ee
Hence large or small values of $M$ for $u$ are exchanged into large or small values of $v(0)$ for $v$ and in the sequel we will essentially express our results using the function $u$.
\subsection{Energy functions}
We consider first the energy function
\bel{Xy-1}\BA {ccc}\displaystyle
r\mapsto H(r)=\myfrac{u^{p+1}}{p+1}+\myfrac{u'^{2}}{2}.
\EA\ee
Then
$$H'(r)=M\abs{u'}^{q+1}-\myfrac{N-1}{r}u'^{2}.
$$
Hence, if $M\leq 0$, $H$ is decreasing, a property often used in \cite{SeZo}. This implies in particular that a radial ground state satisfies
\bel{Xy-2}\BA {ccc}\displaystyle
\abs {u'(r)}\leq \sqrt{\myfrac{2}{p+1}}\left(u(0)\right)^{\frac{p+1}{2}}.
\EA\ee
A similar estimate holds in all the cases.
\bprop{funch} Let $M>0$, $p,q>1$. If $u$ is a radial ground state solution of $(\ref{I-0})$, then the function $H$ defined in $(\ref{Xy-1})$ is decreasing and in particular $(\ref{Xy-2})$ holds.
\es
\Proof Let $u$ be such a radial ground state. By \rprop{zero} we must have $q>\frac N{N-1}$ and
$$\myfrac{r}{u'^{2}}H'=Mr\abs{u'}^{q-1}+1-N\leq \myfrac{(N-1)q-N}{q-1}+1-N=-\frac{1}{q-1},
$$
this implies the claim.\qeda

\subsubsection{Exponential perturbations}
As we have seen it in the introduction, if $q<\frac{2p}{p+1}$ equation $(\ref{I-0})$ can be seen as a perturbation of the Lane-Emden equation $(\ref{I-1})$ while if $q>\frac{2p}{p+1}$ it can be seen as a perturbation of the Ricatti equation $(\ref{I-8})$. Two types of transformations can emphasize these aspects.\\

\nind 1) For $p>1$ set
\bel{X-5}\BA {lll}\displaystyle
u(r)=r^{-\frac{2}{p-1}}x(t),\quad u'(r)=-r^{-\frac{p+1}{p-1}}y(t),\quad t=\ln r,
\EA\ee
then
\bel{X-6}\BA {lll}\displaystyle
x_t=\frac{2}{p-1}x-y\\[2mm]
y_t=-Ky+x^p+Me^{-\gw t}y^q
\EA\ee
with
\bel{X-6a}
K=\frac{(N-2)p-N}{p-1},
\ee
and
\bel{X-7}\BA {lll}\displaystyle
\gw=\myfrac{(p+1)q-2p}{p+1}.
\EA\ee
If $q>\frac{2p}{p+1}$ (resp. $q<\frac{2p}{p+1}$), then $\gw>0$ (resp. $\gw<0$) system $(\ref{X-7})$ is a perturbation of
the Lane-Emden system
\bel{X-8}\BA {lll}\displaystyle
x_t=\frac{2}{p-1}x-y\\[4mm]\displaystyle
y_t=-Ky+x^p,
\EA\ee
at $\infty$ (resp. $-\infty$). The following energy type function introduced in \cite{Lei} is natural with $(\ref{X-8})$
\bel{X-8'}\BA {lll}\displaystyle
\CN(t)=\CL(x(t),y(t))=\myfrac{K}{p-1}x^2-\myfrac{x^{p+1}}{p+1}-\left(\myfrac{2}{p-1}\right)^qMe^{-\gw t}\myfrac{x^{q+1}}{q+1}
-\myfrac{1}{2}\left(\myfrac{2x}{p-1}-y\right)^2,
\EA\ee
and it satisfies
\bel{X-8''}\BA {lll}\displaystyle
\CN'(t)=\left(\myfrac{2x}{p-1}-y\right)\left[L\left(\myfrac{2x}{p-1}-y\right)-Me^{-\gw t}\left(\left(\myfrac{2x}{p-1}\right)^q-y^q\right) \right]\\[4mm]\phantom{-----------------}
+\gw\left(\myfrac{2}{p-1}\right)^qMe^{-\gw t}\myfrac{x^{q+1}}{q+1},
\EA\ee
where $L=N-2-\myfrac{4}{p-1}=K-\myfrac{2}{p-1}$. Relation $(\ref{X-8''})$ will be used later on.
\medskip

\nind 2) For $p,q>1$ set
\bel{X-9}\BA {lll}\displaystyle
u(r)=r^{-\frac{2-q}{q-1}}\xi(t),\quad u'(r)=-r^{-\frac{1}{q-1}}\eta(t),\quad t=\ln r,
\EA\ee
then
\bel{X-10}\BA {lll}\displaystyle
\xi_t=\frac{2-q}{q-1}\xi-\eta\\[4mm]\displaystyle
\eta_t=-\frac{(N-1)q-N}{q-1}\eta+e^{\overline\omega t}\xi^p+M\eta^q
\EA\ee
where
\bel{X-10'}\BA {lll}\displaystyle
\overline\gw=\myfrac{p-1}{q-1}\gw.
\EA\ee
Note that if $q<\frac{2p}{p+1}$ this system at $\infty$ endows the form
\bel{X-11}\BA {lll}\displaystyle
\xi_t=\frac{2-q}{q-1}\xi-\eta\\[4mm]\displaystyle
\eta_t=-\frac{(N-1)q-N}{q-1}\eta+M\eta^q.
\EA\ee
It is therefore autonomous and much easier to study.

\subsubsection{Pohozaev-Pucci-Serrin type functions}
Let $\ga,\gg,\gth,\gk$ be real parameters with $\ga,\gk>0$. Set
\bel{X-12}\BA {lll}\displaystyle
\CZ(r)=r^\gk\left(\myfrac{u'^2}{2}+\myfrac{u^{p+1}}{p+1}+\ga\myfrac{uu'}{r}-\gg u'\abs{u'}^q\right).
\EA\ee
This type of function has been introduced in \cite{SeZo} in their study of equation $(\ref{I-0})$ with $M=1$ with specific parameters. We use it here to embrace all the values of $M$.
We define $\CU$ by the identity
\bel{X-13}\BA {lll}\displaystyle
\CZ'+\gth\abs{u'}^{q-1}\CZ=r^{\gk-1}\CU.
\EA\ee
Then
\bel{X-14}\BA {lll}\displaystyle
\CU=\left(\frac\gk2+\ga+1-N\right)u'^2+\left(\frac\gk{p+1}-\ga\right)u^{p+1}+\ga(\gk-N)\myfrac{uu'}{r}+\left(\myfrac{\gth}{p+1}-\gg q\right)ru^{p+1}\abs{u'}^{q-1}\\[4mm]\displaystyle
\phantom{\CU}
+\left(M+\gg+\frac\gth2\right)r\abs{u'}^{q+1}+\left(\left((N-1)q-\gk\right)^{\phantom{t^;}}\!\!\!\gg-\ga(\gth+M)\right)u\abs{u'}^{q}-\gg(\gth+qM)ru\abs{u'}^{2q-1}.
\EA\ee
\subsection{Some known results in the case $M<0$}
We recall the results of \cite{ChWe}, \cite{SeZo} and \cite{PoQuSo} relative to the case $M<0$.
\bth{known} 1- Let $N\geq 3$ and $1<p\leq \frac{N}{N-2}$.\\
1-(i) If $q>\frac{2p}{p+1}$, there is no ground state  for any $M<0$ (\cite[Theorem C]{SeZo}).\\
1-(ii) If $1<q<\frac{2p}{p+1}$ there exists a ground state when $|M|$ is large \cite[Proposition 5.7]{ChWe} and there exists no ground state
when $|M|$ is small (\cite{PoQuSo}).\smallskip

\nind 2- Assume $\frac {N}{N-2}<p<\frac {N+2}{N-2}$ and let $\overline q$ be the unique root in $(\frac{2p}{p+1},p)$ of the quadratic equation
$$(N-1)(X-p)^2-(N+2-(N-2)p)((p+1)X-2p)X=0.
$$
2-(i) If $\overline q\leq q<p$ there exists no ground state  for any $M<0$ (\cite[Theorem C]{SeZo}).\\
2-(ii) If $\frac{2p}{p+1}<q<\overline q$, there exists no ground state for $|M|$. It is an open question whether  there could exist a finite number of $M$ for which there exists a ground state (\cite[Theorem 4]{SeZo}). \\
2-(iii) If $1<q<\frac{2p}{p+1}$, there exists a ground state for large $|M|$ (\cite[Proposition 5.7]{ChWe}) and no ground state when $|M|$ is small (\cite{PoQuSo}).\smallskip

\nind 3- Assume $p>\frac{N+2}{N-2}$ and $q>1$ and let $Q_{N,p}=\frac{2(N-1)p}{2N+p+1}\in (\frac{2p}{p+1},p)$.\\
3-(i) If $Q_{N,p}<q<p$ there exists a ground state for $|M|$ small. \\
3-(ii) If $1<q\leq Q_{N,p}$ there exists a ground state for any $M<0$ (\cite[Theorem A]{SeZo}). \smallskip

\nind 4- Assume $p=\frac{N+2}{N-2}$. There exists at least one $M<0$ such that there exists a ground state if and only if $1<q<p$. More precisely:\\
4-(i) If $\frac{2p}{p+1}<q<p$ there exists ground state if $|M|$ is small (\cite[Theorem B]{SeZo}).\\
4-(ii) If $q\geq \frac{2p}{p+1}$ there exists a ground state for any $M<0$ (\cite[Theorem A]{SeZo}).
\es

\nind\Remark It is interesting to quote that when $M<0$ and $q\geq \frac{2p}{p+1}$, there holds \cite[Theorem 3]{SeZo},
$$u(r)=O(r^{-\frac{2}{p-1}})\quad\text{and }\;u'(r)=O(r^{-\frac{p+1}{p-1}})\quad\text{when }r\to\infty.
$$
\subsection{The case $M>0$}
The next result is a consequence of Theorem A.
\bth{surcri} Let $M>0$, $p>1$ and $q>\frac{2p}{p+1}$ then there exists no radial ground state satisfying $u(0)=1$ when $M$ is large.
\es
\Proof Suppose that such a solution $u$ exists. From Theorem A and \rprop{zero} there holds
\bel{X-15}\BA {lll}\displaystyle
\sup_{r>0}\abs{u'(r)}\leq c_{N,p,q}|M|^{-\frac{p+1}{(p+1)q-2p}}\quad\text{and }\; \sup_{r>0}r^{\frac{p+1}{p-1}}\abs{u'(r)}\leq c_{N,p}.
\EA\ee
As a consequence, if $r>R>0$,
$$\BA {lll}
1-u(r)=u(0)-u(r)=u(0)-u(R)+u(R)-u(r)\leq c_{N,p,q}|M|^{-\frac{p+1}{(p+1)q-2p}}R+\myint{R}{\infty}\abs{u'(s)}ds\\[4mm]
\phantom{1-u(r)}
\leq c_{N,p,q}|M|^{-\frac{p+1}{(p+1)q-2p}}R+c'_{N,p}R^{-\frac{2}{p-1}},
\EA$$
with $c'_{N,p}=\frac{p-1}{2}c_{N,p}$. Since $u(r)\to 0$ when $r\to\infty$, we take $R=|M|^{\frac{p-1}{(p+1)q-2p}}$ and derive
\bel{X-16}\BA {lll}\displaystyle
1\leq \left(c_{N,p,q}+c'_{N,p}\right)|M|^{-\frac{2}{(p+1)q-2p}},
\EA\ee
and the conclusion follows.\qeda\medskip

\nind \Remark If we use \rprop{funch} we can make estimate $(\ref{X-16})$ more precise.
\subsubsection{The case $M>0$, $1<p\leq\frac{N+2}{N-2}$}
It is a consequence of our general results that  there is no radial ground state for large $M$ or for small $M$ when $1<q\leq \frac{2p}{p+1}$ and $1<p<\frac{N+2}{N-2}$. Indeed, if $1<q<\frac{2p}{p+1}$ is a consequence of the equivalence statement between a priori estimate and non-existence of ground state proved in \cite{PoQuSo}, and if $q= \frac{2p}{p+1}$ it follows from Theorems C and E. Actually in the radial case, the result is more general.
\bth{well} Let $M>0$ and $1< p<\frac{N+2}{N-2}$. If $1<q\leq p$, there exists no radial ground state for any $M$. If $q>p$ there exists no radial ground state for $M$ small enough.
\es
\Proof By \rprop{zero}, we may assume  $N\geq 3$ and
\bel{X-17}\frac{N}{N-2}<p\leq \frac{N+2}{N-2}\quad\text{and }\;q>\frac{N}{N-1}.\ee
(i) Assume first $q< \frac{2p}{p+1}$. We use the system $(\ref{X-6})$. Then $\gw$, defined by $(\ref{X-7})$ is negative. Hence the Leighton function $\CN$ defined by $(\ref{X-8'})$ is nonincreasing since $L\leq 0$ when $p\leq\frac{N+2}{N-2}$. Furthermore since
$(x(t),y(t))\to (0,0)$ when $t\to-\infty$ and $e^{-\gw t}\to 0$, we get $\CN(-\infty)=0$ it follows that $\CN(t)<0$ for $t\in\BBR$.
Moreover, by \rprop{zero},
$$u(r)=O(r^{-\frac{2-q}{q-1}})\quad\text{as }r\to\infty\Longleftrightarrow x(t)=O(e^{\frac{q(p+1)-2p}{(p-1)(q-1)}t})=o(1)\quad\text{as }t\to\infty
$$
This implies  $e^{-\gw t}x^{q+1}(t)=O(e^{2\frac{q(p+1)-2p}{(p-1)(q-1)}t})=o(1)$ as $t\to\infty$ and $\CN(\infty)=0$, contradiction.\smallskip

\nind (ii) Assume next $\frac{2p}{p+1}\leq q\leq p$. We consider the function $(\ref{X-12})$ with the parameters
$$\gk=\myfrac{2(p+1)(N-1)}{p+3}=(p+1)\ga\quad\text{and }\;\gg=-\myfrac{2M}{q(p+1)+2}=\myfrac{\gth}{q(p+1)},
$$
already used by \cite{SeZo} when $M=-1$, and we get with $\CU$ defined by $(\ref{X-13})$,
$$\CU=\myfrac{2}{(p+3)^2}\myfrac{u\abs{u'}}{r}\left(A+BM\chi+CM\chi^2\right)\quad\text{with}\;\chi=\myfrac{p+3}{2+q(p+1)}r\abs{u'}^{q-1},
$$
where
\bel{X-18}
A=(N-1)(N+2-(N-2)p)\,,\; B=2(N-1)(p-q)\,,\; C=q(q(p+1)-2p).
\ee
By our assumptions $A\geq 0$, $B\geq 0$ and $C>0$. Hence $\CU>0$. This implies
$$\CZ(r) =e^{-\int_{0}^{r}\gth\abs{u'}^{q-1}ds}\CZ(0)+\myint{0}{r}e^{-\gth\int_{s}^{r}\abs{u'}^{q-1}d\gs}s^{\gk-1}\CU(s)ds
=\myint{0}{r}e^{-\gth\int_{s}^{r}\abs{u'}^{q-1}d\gs}s^{\gk-1}\CU(s)ds,
$$
since $\CZ(0)=0$. If $u$ is a ground state, then $u'(r)\to 0$ as $r\to\infty$, thus $u\abs{u'}^q\leq u\abs{u'}^{\frac{2p}{p+1}}$. Hence, from \rprop{zero}, $u'^2(r)=O(r^{-2\frac{p+1}{p-1}})$ as $r\to\infty$. The other terms $u^{p+1}(r)$, $r^{-1}u(r)u'(r)$ and $u\abs{u'}^{\frac{2p}{p+1}}$ satisfy the same bound, hence
$$\CZ(r)=O(r^{\gk-\frac{2(p+1)}{p-1}})=O(r^{\frac{2(p+3)(N-1)}{p+3}-\frac{2(p+1)}{p-1}})=O(r^{\frac{2(p+1)((N-2)p-(N+2))}{(p+3)(p-1)}}).
$$
Then $\CZ(r)\to 0$ when $r\to\infty$, contradiction. \smallskip

\nind (iii) Suppose $q>p$ and $u$ is a ground state. By \rprop{funch} and $(\ref{X-15})$, there holds
$$r\abs{u'}^{q-1}=r\abs{u'}^{\frac{p-1}{p+1}}\abs{u'}^{q-\frac{2p}{p+1}}\leq c_{N,p}.
$$
Then $\chi=\frac{p+3}{2+q(p+1)}r\abs{u'}^{q-1}\leq c_{N,p}$. Hence, if $M\leq M_{N,p}$ for some $M_{N,p}>0$, $\CU$ is positive as
$A$ is. We conclude as above. \qeda\medskip


\subsubsection{The case $M>0$, $p>\frac{N+2}{N-2}$}

We  recall that in Theorem C if $q=\frac{2p}{p+1}$ and $p>1$ there is no ground state whenever $M>M_{N,p}$, see $(\ref{I-14b})$. In  Theorem A' if $1<q<\frac{2p}{p+1}$ and $p>1$ there is no ground state $u$ such that $u(0)=1$ if $M$ is too large.
In the next result we complement \rth{surcri} for small value of $M$ in assuming $q>\frac{2p}{p+1}$.

\bth{surcri2} If $p>\frac{N+2}{N-2}$ and $q\geq\frac{2p}{p+1}$ then there exist radial ground states for $M>0$ small enough.
\es
\Proof First we consider the function $\CZ$ with $k=N$ and obtain
$$\CZ(r)=r^{N}\left(\myfrac{u'^2}{2}+\myfrac{u^{p+1}}{p+1}+\ga\myfrac{uu'}{r}-\gg u\abs{u'}^q\right).
$$
The function vanishes at the origin. We compute $\CU$ from the identity $\CZ'+\gth\abs{u'}^{q-1}\CZ=r^{N-1}\CU$ and get
$$\BA {lll}
\CU=\left(\ga-\myfrac{N-2}{2}\right)u'^2+\left(\myfrac{N}{p+1}-\ga\right)u^{p+1}+\left(\myfrac{\gth}{p+1}-\gg q\right)ru^{p+1}\abs{u'}^{q-1}\\[4mm]
\phantom{--,}
+\left(M+\gg+\myfrac{\gth}{2}\right)r\abs{u'}^{q+1}+\left[\left((N-1)q-N\right)\gg-\ga(\gth +M)^{\phantom{\frac{A^a}{B‰}}}\!\!\!\!\!\!\!\right]u\abs{u'}^{q}
-\gg(\gth+qM)ru\abs{u'}^{2q-1}.
\EA$$
If $\gg=0$ and $\gth=-2M$, then
$$\CU=\left(\ga-\myfrac{N-2}{2}\right)u'^2+\left(\myfrac{N}{p+1}-\ga\right)u^{p+1}-\myfrac{2M}{p+1}ru^{p+1}\abs{u'}^{q-1}
+\ga Mu\abs{u'}^{q}.
$$
If $u$ is a regular solution which vanishes at some $r_0>0$, then $\CZ(r_0)=2^{-1}r_0^2u'^N(r_0)>0$. As $p>\frac{N+2}{N-2}$, by choosing $\ga=\frac{1}{2}\left(\frac{N}{p+1}+\frac{N-2}{2}\right)$ we have $\frac{N}{p+1}<\ga<\frac{N-2}{2}$. We define $\ell>0$ by  $(N-2)p-(N+2)=4(p+1)\ell$, then
$\frac{N-2}{2}-\ga=\ga-\frac{N}{p+1}=\ell$ and then
$$\CU\leq -\ell (u'^2+u^{p+1})+M\ga u\abs{u'}^q.
$$

Assume first $q<2$, we have from H\"older's inequality and $0<r\leq r_0$ where $u$ is positive
$$u\abs{u'}^q\leq \myfrac{q}{2}u'^2+\myfrac{2-q}{2}\abs{u}^{\frac{2}{2-q}}\leq u'^2+\abs{u}^{\frac{2}{2-q}},
$$
and
$$\BA {lll}
\CU+(\ell-M)u'^2\leq M\ga u^{\frac{2}{2-q}}-\ell u^{p+1}=\ell u^{p+1}\left(\myfrac{M\ga}{\ell}u^{\frac{q(p+1)-2p}{2-q}}-1\right)
\leq \ell u^{p+1}\left(\myfrac{M\ga}{\ell}-1\right)
\EA$$
since $q\geq \frac{2p}{p+1}$ and $u\leq u(0)=1$. Taking $M\leq\frac{\ell}{\ga}=\frac{(N-2)p-N-2}{(N-2)p+3N-2}$, $\CU$ is negative implying that  $r\mapsto e^{-2M\int_0^r|u'|^{q-1}ds}\CZ(r)$ is nonincreasing. Since $\CZ(0)=0$, $\CZ(r)\leq 0$, a  contradiction.

If $q=2$, then $\CU\leq -\ell (u'^2+u^{p+1})+M\ga u'^2$ since $u\leq 1$ on $[0,r_0]$. We still infer that $\CU\leq 0$ if we choose   $M\leq\frac{\ell}{\ga}$.

Finally, if $q>2$, we have from Theorem A, $u'\leq C_{N,p,q} M^{-\frac{p+1}{(p+1)q-2p}}$. Therefore, using again the decay of $u$ from
$u(0)=1$,
$$M\ga u\abs{u'}^q\leq M\ga u\abs{u'}^{q-2}u'^2\leq M\ga C^{q-2}_{N,p,q} M^{-\frac{(p+1)(q-2)}{(p+1)q-2p}}u'^2=
\ga C^{q-2}_{N,p,q}M^{\frac{2}{(p+1)q-2p}}u'^2.
$$
Hence $\CU\leq -\left(\ell-\ga C^{q-2}_{N,p,q}M^{\frac{2}{(p+1)q-2p}}\right)u'^2$. Taking
$$M^{\frac{2}{(p+1)q-2p}}\leq C^{2-q}_{N,p,q}\frac{(N-2)p-N-2}{(N-2)p+3N-2}$$
we conclude that $\CU<0$ which ends the proof as in the previous cases. \qeda\medskip

Theorem F is the combination of \rth{surcri}, \rth{well} and \rth{surcri2}.

\mysection{Separable solutions}

We denote by $(r,\gs)\in\BBR_+\ti S^{N-1}$ the spherical coordinates in $\BBR^N$. Then equation $(\ref{I-0})$
takes the form
\bel{Sep0}
\BA{lll}
- u_{rr}    -\myfrac{N-1}{r} u_{r}   -\myfrac{1}{r^2}\Gd'u=|u|^{p-1}+M\left(u_{r} ^2+\myfrac{1}{r^2}\abs{\nabla'u}^2\right)^{\frac{q}{2}},
\EA
\ee
where $\Gd'$ is the Laplace-Beltrami operator on $S^{N-1}$ and $\nabla'$ the tangential gradient.
If we look for separable nonnegative solutions of $(\ref{I-0})$ i.e. solutions under the form $u(r,\gs)=\psi(r)\gw(\gs)$, then $q=\frac{2p}{p+1}$,  $\psi(r)=r^{-\frac{2}{p-1}}$, and $\gw$ is a solution of
\bel{Sep1}
\BA{lll}
-\Gd'\gw+\myfrac{2K}{p-1}\gw=\gw^p+M\left(\left(\myfrac{2}{p-1}\right)^2\gw^2+\abs{\nabla'\gw}^2\right)^{\frac{p}{p+1}},
\EA
\ee
where
$K$ is defined in $(\ref{X-6a})$.
Throughout this section we assume
\bel{Sep2'}
p>1\quad\text{and }\; q=\frac{2p}{p+1}.        \ee

\subsection{Constant solutions}
The constant function $\gw=X$ is a solution of \eqref{Sep1} if
\bel{Sep3}
\BA{lll}
X^{p-1}+M\left(\myfrac{2}{p-1}\right)^{\frac{2p}{p+1}}X^{\frac{p-1}{p+1}}-\myfrac{2K}{p-1}=0.
\EA
\ee
For $N=1, 2$ and $p>1$ or $N\geq 3$ and $1<p<\frac{N}{N-2}$, we recall that $\gm^*=\gm^*(N)$ has been defined in $(\ref{I-14})$.
The following result is easy to prove
\bprop{const1} 1- Let $M\geq 0$ then there exists a unique positive root $X_M$ to $(\ref{Sep3})$ if and only if $p>\frac{N}{N-2}$. Moreover the mapping
$M\mapsto X_M$ is continuous and decreasing from $[0,\infty)$ onto $(0,\left(\frac{2K}{p-1}\right)^{\frac{1}{p-1}}]$.\smallskip

\nind 2- Let $M<0$, $N\geq 3$ and $p\geq\frac{N}{N-2}$ then there exists a unique positive root $X_M$ to $(\ref {Sep3})$ and the mapping
$M\mapsto X_M$ is continuous and decreasing from $(-\infty,0]$ onto $[\left(\frac{2K}{p-1}\right)^{\frac{1}{p-1}},\infty)$.\smallskip

\nind 3- Let $M<0$, $N=1,2$ and $p>1$ or $N\geq 3$ and $1<p<\frac{N}{N-2}$ then there exists no positive root to $(\ref {Sep3})$ if
$-\gm^*<M\leq 0$. If $M=M^*:=-\gm^*$ there exists a unique positive root $X_{M^*}=\left(\frac{2|K|}{p(p-1)}\right)^{\frac{1}{p-1}}$. If $M<-\gm^*$ there exist two positive roots
$X_{1,M}<X_{2,M}$. The mapping $M\mapsto X_{1,M}$ is continuous and increasing from $(-\infty,\gm^*]$ onto
$(0,X_{M^*}]$.  The mapping $M\mapsto X_{2,M}$ is continuous and decreasing from $(-\infty,\gm^*]$ onto
$[X_{M^*},\infty)$.
\es
 \nind{\it Abridged proof}. Set
\bel{x1}
f_M(X)=X^{p-1}+M\left(\myfrac{2}{p-1}\right)^{\frac{2p}{p+1}}X^{\frac{p-1}{p+1}}-\myfrac{2K}{p-1},
\ee
then $f'_M(X)=(p-1)X^{p-2}+M\frac{p-1}{p+1}\left(\myfrac{2}{p-1}\right)^{\frac{2p}{p+1}}X^{-\frac{2}{p+1}}$. \smallskip

\nind 1- If $M$ is nonnegative, $f_M$ is increasing from $-\frac{2K}{p-1}=-\frac{2[(N-2)p-N]}{(p-1)^2}$ to $\infty$; hence, if $p>\frac{N}{N-2}$ there exists a unique $X_M> 0$ such that $f_M(X_M)=0$, while if $1<p<\frac{N}{N-2}$, $f_M$ admits no zero on $[0,\infty)$. Since $f_M>f_{M'}$ for $M>M'>0$, there holds $X_M>X_{M'}$, By the implicit function theorem the mapping $M\mapsto X_M$ is $C^1$ and decreasing from $[0,\infty)$ onto $(0,\left(\frac{2K}{p-1}\right)^{\frac{1}{p-1}}]$. Actually it can be proved that  (see \cite[Proposition 2.2]{BiGaVe2})
\bel{x31}
X_M=\frac{p-1}{2}\left(\frac KM\right)^{\frac{p+1}{p-1}}(1+o(1))\quad\text{as }\;M\to \infty.
\ee

\nind 2- If $M$ is negative, $f_M$ achieves it minimum on $[0,\infty)$ at $X_0=\left(\frac{-M}{p+1}\right)^\frac{p+1}{p(p-1)}\left(\frac{2}{p-1}\right)^\frac{2}{p-1}$, and
$$\BA {lll}f_M(X_0)=-\myfrac{p}{\left(p+1\right)^{\frac {p+1}p}}\left(\myfrac{2}{p-1}\right)^{2}(-M)^{\frac{p+1}{p}}-\myfrac{2K}{p-1}\\[4mm]
\phantom{f_M(X_0)}
=-\left(\myfrac{2}{p-1}\right)^{2}\left(\myfrac{p}{\left(p+1\right)^{\frac {p+1}p}}(-M)^{\frac{p+1}{p}}+\myfrac{(N-2)p-N}{2}\right).
\EA$$
Since $K>0$, there exists a unique $X_M>0$ such that $f_M(X_M)=0$ and $X_M>X_0$. The mapping $M\mapsto X_M$ is $C^1$ and decreasing from $(-\infty, 0]$ onto $[\left(\frac{2K}{p-1}\right)^{\frac{1}{p-1}},\infty)$. The following estimate holds
\bel{x35}\BA {lll}\displaystyle
\max\left\{\left(\frac{2K}{p-1}\right)^{\frac{1}{p-1}},\left(\frac{2}{p-1}\right)^{\frac{2}{p-1}}|M|^{\frac{p+1}{p(p-1)}}\right\}\leq X_M\\[4mm]
\displaystyle\phantom{---------}\leq 2^{\frac{2}{p-1}}\left(\left(\frac{2K}{p-1}\right)^{\frac{1}{p-1}}+\left(\frac{2}{p-1}\right)^{\frac{2}{p-1}}|M|^{\frac{p+1}{p(p-1)}}\right).
\EA\ee
\smallskip

\nind 3- If $N=1,2$ and $p>1$ or $N\geq 3$ and $1<p<\frac{N}{N-2}$, then $f_M(0)>0$. Hence, if $f_M(X_0)>0$ there exists no positive
root to $f_M(X)=0$. Equivalently, if $-\gm^*<M<0$. If $f_M(X_0)=0$, $X_0$ is a double root and this is possible only if $M=-\gm^*$, hence $X_{-\gm^*}=\left(\frac{2|K|}{p(p-1)}\right)^{\frac{1}{p-1}}$. If $f_M(X)<0$, or equivalently, if $M<-\gm^*$, the equation $f_M(X)=0$ admits two positive roots $X_{1,M}<X_0<X_{2,M}$. The monotonicity of the $X_{j,M}$, j=1,2, and their range follows easily from the monotonicity of $M\mapsto f_M(X)$ for $M<0$.  Actually the following asymptotics hold when $M\to -\infty$,
\bel{x33}
X_{1,M}=\frac{p-1}{2}\left(\frac{K}{M}\right)^{\frac{p+1}{p-1}}(1+o(1))\;\text{and }\;\, X_{2,M}=\left(\frac{2}{p-1}\right)^\frac{2}{p-1}\left(-M\right)^{\frac{p+1}{p(p-1)}}(1+o(1)).
\ee
\qeda\medskip


 \subsection{Bifurcations}
 We set
                   \bel{Sep4}
            \BA{lll}
A(\gw)=   -\Gd'\gw+\myfrac{2K}{p-1}\gw-\gw^p-M\left(\left(\myfrac{2}{p-1}\right)^2\gw^2+\abs{\nabla'\gw}^2\right)^{\frac{p}{p+1}},
                   \EA
                   \ee
   If $\eta\in C^\infty(S^{N-1})$ and if there exists a constant positive solution $X$ to $A(X)=0$ we have
   $$     \BA{lll}
   \myfrac{d}{d\gt}A(X+\gt\eta)\lfloor_{\gt=0}=-\Gd'\eta+\left(\myfrac{2K}{p-1}-pX^{p-1}-M\myfrac{2p}{p+1}\left(\myfrac{2}{p-1}\right)^{\frac{2p}{p+1}}X^{\frac{p-1}{p+1}}\right)\eta.
   \EA
   $$
   Hence the problem is singular if
                    \bel{Sep5}
            \BA{lll}
-\myfrac{2K}{p-1}+pX^{p-1}+M\myfrac{2p}{p+1}\left(\myfrac{2}{p-1}\right)^{\frac{2p}{p+1}}X^{\frac{p-1}{p+1}}=\gl_k,
                   \EA
                   \ee
                   where $\gl_k=k(k+N-2)$ is the k-th nonzero eigenvalue of $-\Gd'$ in $H^1(S^{N-1})$.    The following result follows classically from the standard bifurcation theorem from a simple eigenvalue  (which can always be assumed if we consider functions   depending only on the azimuthal angle on $S^{N-1}$ reducing the eigenvalue problem to a simple Legendre type ordinary differential equation) see e.g. \cite[Chapter 13]{Smo} and identity $(\ref{Sep3})$.

 \bth{bif} Let $M_0\in\BBR$ and $X_{M_0}$ be a  constant solution of $(\ref{Sep1})$. If $X_{M_0}$ satisfies for some $k\in\BBN^*$,
                       \bel{Sep7}
{M_0}\left(\myfrac{2}{p-1}\right)^{\frac{2p}{p+1}}X_{M_0}^{\frac{p-1}{p+1}}=\myfrac{p+1}{p(p-1)}\left(2K-\gl_k\right),
                   \ee
there exists a continuous branch of nonconstant positive solutions $(M,\gw_M)$ of $(\ref{Sep1})$ bifurcating from the $({M_0},X_{M_0})$.
 \es

Since  $M\left(\myfrac{2}{p-1}\right)^{\frac{2p}{p+1}}X_M^{\frac{p-1}{p+1}}=\myfrac{2K}{p-1}-X_M^{p-1}$ by $(\ref{Sep3})$
         the following statements follow immediately from \rprop{const1}.
   \blemma{var} Set $\Gf(M)=M\left(\frac{2}{p-1}\right)^{\frac{2p}{p+1}}X_M^{\frac{p-1}{p+1}}$ when  $X_M$ is uniquely determined, and
$\Gf_j(M)=M\left(\frac{2}{p-1}\right)^{\frac{2p}{p+1}}X_{j,M}^{\frac{p-1}{p+1}}$, j=1,2, if there exist two equilibria. Then   \smallskip

   \nind 1- If $N\geq 3$ and $p>\frac{N}{N-2}$, the mapping $M\mapsto \Gf(M)$ is continuous and increasing from $[0,\infty)$ onto $[0,\frac{2K}{p-1})$.\smallskip

   \nind 2- If $N\geq 3$ and $p\geq\frac{N}{N-2}$, the mapping $M\mapsto \Gf(M)$ is continuous and increasing from $(-\infty,0]$ onto $(-\infty,0]$.\smallskip

   \nind 3- If $N=1,2$ and $p>1$ or $N\geq 3$ and $1<p<\frac{N}{N-2}$, the mapping $M\mapsto \Gf_1(M)$ (resp $M\mapsto \Gf_2(M)$) is continuous and decreasing (resp. increasing) from $(-\infty,-\gm^*]$ onto $[\frac{2K}{p-1},0)$ (resp. $(-\infty,\frac{2K}{p-1}]$).
 \es

 If we analyse the range $R[\Gf]$ of $\Gf$ or $R[\Gf_j]$ of $\Gf_j$, we prove the following  result.

 \bth{bifth} 1- Let  $N\geq 3$ and $p\geq \frac{N}{N-2}$. \smallskip

 \nind 1-(i) There exists a continuous curve of bifurcation $(M,\gw_M)$ issued from  $(M_0,X_{M_0})$ for some $M_0=M_0(p)\geq 0$ if and only if $p\geq\frac{N+1}{N-3}$ and $k=1$.  Furthermore $M_0(\frac{N+1}{N-3})=0$. \smallskip

 \nind 1-(ii) The bifurcation curve $s\mapsto (M(s),\gw_{M(s)})$, is defined on $(-\ge_0,\ge_0)$ for some $\ge_0>0$ and verifies  $(M(0),\gw_{M(0)})=(M_0,X_{M_0})$.\smallskip

   \nind 2- Let $N\geq 3$ and $p\geq\frac{N}{N-2}$. \smallskip

 \nind 2-(i) For any $k\geq 1$ there exist $M_k<0$ and a continuous branch of bifurcation $(M,\gw_M)$ issued from  $(M_k,X_{M_k})$,   with the restriction that $p< \frac{N+1}{N-3}$ if $k=1$. \smallskip

 \nind 2-(ii) The bifurcation curve $s\mapsto (M(s),\gw_{M(s)})$, is defined on $(-\ge_0,\ge_0)$ for some $\ge_0>0$ and verifies  $(M(0),\gw_{M(0)})=(M_0,X_{M_0})$. Finally $M_k\to-\infty$ when $k\to\infty$.\smallskip

   \nind 3- let $N=1,2$ and $p>1$, or $N\geq 3$ and $1<p<\frac{N}{N-2}$. \smallskip

 \nind 3-(i) There exists no $M<0$ such that $\frac{2K}{p-1}<\Gf_1 (M)<0$, and a countable set of $M_k<0$, $k\geq 1$, such that $\Gf_2(M_k)=\frac{p+1}{p(p-1)}\left(2K-\gl_k\right)$. \smallskip

 \nind 3-(ii) There exist a countable branches of bifurcation of solutions $(M_k(s),\gw_{M_k(s)})$ issued from $(M_k, X_{2,M_k})$.
 \es
 \Proof  {\it Assertion  1}. Since from \rlemma {var}, $R[\Gf]=[0,\frac{2K}{p-1})$ for $M\geq 0$, we have to see under what condition on $p\geq \frac{N}{N-2}$ one can find  $k\geq 1$ such that
 $$0\leq \myfrac{p+1}{p(p-1)}\left(2K-\gl_k\right)<\frac{2K}{p-1}\Longleftrightarrow \frac {2K}{p+1} <\gl_k\leq 2K.
 $$
 Since $K<N$ and $\gl_k\geq 2N$ for $k\geq 2$, the only possibility for this last inequality to hold is $k=1$. The inequality $\frac {2K}{p+1} <N-1$ always holds since $p>1$, while the inequality $N-1=\gl_1\leq 2K$ is equivalent to $p\geq \frac{N+1}{N-3}$. Therefore $M_0=0$ and $X_{M_0}=\left(\frac{2K}{p-1}\right)^{\frac 1{p-1}}$. If we consider only functions on the sphere $S^{N-1}$  which depend uniquely on the azimuthal angle $\gth=\tan^{-1}(x_N\lfloor_{S^{N-1}})$, the function $\psi_1(\gs)= x_N\lfloor_{S^{N-1}}$ is a eigenfunction of $-\Gd'$ in
 $H^1(S^{N-1})$ with multiplicity one. Hence the bifurcation branch is locally a regular curve $s\mapsto (M(s),\gw_{M(s)})$ with  $0\leq s<\ge'_0$ and we construct the bifurcating solution on  $S^{N-1}$ using the classical tangency condition \cite[Theorem 13.5]{Smo},
\bel{Sep7a}\gw_{M(s)}=X_{M_0}+s(\psi_1+\gz_s)
\ee
where $\gz_s\in H^{1}(S^{N-1})$, is orthogonal to $\psi_1$ in $H^{1}(S^{N-1})$ and  satisfies
$\norm{\gz_s}_{C^1}=o(1)$ when $s\to 0$.
 This implies the claim.
\smallskip

 \nind {\it Assertion  2}. Since $R[\Gf]=(-\infty,0)$ for $M< 0$, we have to find $k\geq 1$ such that
 $$\myfrac{p+1}{p(p-1)}\left(2K-\gl_k\right)< 0\Longleftrightarrow 2K< \gl_k.
 $$
 As in Case 1, $K<2N$, then inequality $2K\leq\gl_k$ holds for all $k\geq 2$, and if $k=1$ this is possible only if $p< \frac{N+1}{N-3}$.
 The construction of the bifurcating curve is the same as in Case 1.
 \smallskip

 \nind {\it Assertion   3}. We have $R[\Gf_1]=[\frac{2K}{p-1},0)$ for $M\leq-\gm^*$. If we look for the existence of some
 $k\geq 1$ such that
 $$\frac{2K}{p-1}\leq \myfrac{p+1}{p(p-1)}\left(2K-\gl_k\right)<0\Longleftrightarrow 2K\leq \gl_k<\frac{2K}{p+1};
 $$
we get an impossibility since $K<0$. Hence there exists no $M_0<0$ such that $(M_0,X_{1,M_0})$ is a bifurcation point. We have also
$R[\Gf_2]=(-\infty,\frac{2K}{p-1}]$ for $M\leq-\gm^*$. Now the condition for the existence of a bifurcation branch issued from $(M_0,X_{2,M_0})$ for some
$M_0\leq -\gm^*$ is
 $$\myfrac{p+1}{p(p-1)}\left(2K-\gl_k\right)\leq \frac{2K}{p-1}\Longleftrightarrow \gl_k\geq \frac{2K}{p+1},
 $$
 which is always true for any $k\geq 1$ and $1<p<\frac N{N-2}$.\qeda\medskip

 \nind \Remark The exponent $p=\frac{N+1}{N-3}$ is the Sobolev critical exponent on $S^{N-1}$. If one consider the Lane-Emden equation
 with a Leray potential
 \bel{EF}
 -\Gd u+\gl |x|^{-2}u=u^{\frac{N+1}{N-3}},
 \ee
with $\gl\in\BBR$, then the separable solutions $u(r,\gs)=r^{-\frac{N-3}{2}}\gw(\gs)$ verify
 \bel{EFb}
-\Gd'\gw+\left(\frac{(N-1)(N-3)}{4}-\gl\right)\gw-\gw^{\frac{N+1}{N-3}}=0\quad\text{on }\; S^{N-1}.
 \ee
 It was observed in \cite{BiVe} that there exists a branch of bifurcation $(\gl,\gw_\gl)$ with $\gl>0$ issued from $(0,\gw_0)$, where
$ \gw_0$ is the constant explicit solution of $( \ref{EFb})$.\medskip

\nind\Remark In \rth{bifth}-1- and the above remark, we conjectured that on the bifurcating curve there holds locally $M(s)<M_0$, and that for any
$p\geq \frac{N+1}{N-3}$  there exists $M_0:=M_0(p)$ such that for $M> M_0$ all the positive solutions to $(\ref{Sep1})$ are constant, furthermore $M_0$ is defined by $(\ref{Sep7})$. When $p=\frac{N+1}{N-3}$, then $M=0$ and there exists infinitely many positive solutions to $(\ref{Sep1})$ \cite[Proposition 5.1]{BiVe}. When $\frac{N}{N-2}<p< \frac{N+1}{N-3}$, it is unclear if the branches of bifurcation $(M(s),\gw_{M(s)})$ issued from $(M_0,\gw_{M_0})$ with $M_0<0$ are such that
$M(s)$ keeps a constant sign. If it is the case one could expect that if $M\geq 0$ and $\frac{N}{N-2}<p< \frac{N+1}{N-3}$, all the positive solutions to $(\ref{Sep1})$ are constant.\medskip

The following statement is an immediate consequence of \rth{bifth}.

\bcor{bifcior} 1-If $p>1$ and $q=\frac{2p}{p+1}$ there always exist nonradial positive singular solutions of $(\ref{I-0})$ in $\BBR^N\setminus\{0\}$ under the form $u(r,\gs)=r^{-\frac{2}{p-1}}\gw(\gs)$.

\nind 2- If $N\geq 4$ and $p>\frac{N+1}{N-3}$, the functions are obtained by bifurcation from $X_M$ with $M>0$.

\nind 3- If $N\geq 3$ and $\frac{N}{N-2}\leq p<\frac{N+1}{N-3}$, the functions are obtained by bifurcation from $X_M$ with $M<0$.

\nind 4-If $N=1,2$ and $p>1$ or $N\geq 3$ and $1<p<\frac{N}{N-2}$, the functions are obtained by bifurcation from $(M_k, X_{2,M_k})$ with $M_k<-\gm^*$ and $k\geq 1$.
\es

\end{document}